\documentclass[10pt,sumlimits,namelimits]{article}

\usepackage{amsthm}
\usepackage{amsmath}
\usepackage{amsfonts}
\usepackage{amssymb}
\usepackage{mathrsfs}
\usepackage{enumerate}
\usepackage{color}
\usepackage{graphicx}
\usepackage[colorlinks=true,urlcolor=blue]{hyperref}
\usepackage[applemac]{inputenc}

\newcommand{\A}{\mathbb{A}}
\newcommand{\B}{\mathbb{B}}
\newcommand{\N}{\mathbb{N}}
\newcommand{\Z}{\mathbb{Z}}

\newcommand{\R}{\mathbb{R}}
\newcommand{\C}{\mathbb{C}}

\renewcommand{\S}{\mathbb{S}}
\newcommand{\vA}{{\cal A}}
\newcommand{\vB}{{\cal B}}
\newcommand{\vC}{{\cal C}}

\newcommand{\vO}{{\cal O}}

\newcommand{\Br}{\mathscr{B}}
\newcommand{\Dr}{\mathscr{D}}
\newcommand{\0}{\bs{0}}
\newcommand{\va}{\bs{a}}
\newcommand{\vb}{\bs{b}}
\newcommand{\f}{\bs{f}}
\newcommand{\g}{\bs{g}}
\newcommand{\F}{\bs{F}}
\newcommand{\G}{\bs{G}}

\renewcommand{\u}{\bs{u}}
\renewcommand{\v}{\bs{v}}
\newcommand{\vw}{\bs{w}}
\newcommand{\z}{\bs{z}}
\newcommand{\vphi}{\varphi}
\newcommand{\eps}{\varepsilon}

\newcommand{\dsp}{\displaystyle}
\newcommand{\ovl}{\overline}
\newcommand{\udl}{\underline}

\newcommand{\vmax}{\max\limits}
\newcommand{\vmin}{\min\limits}
\newcommand{\vsup}{\sup\limits}

\newcommand{\vint}{\int\limits}
\newcommand{\vsum}{\sum\limits}
\newcommand{\inj}{\hookrightarrow}
\newcommand{\tends}{\longrightarrow}
\newcommand{\weak}{\rightharpoonup}
\newcommand{\wt}{\widetilde}

\newcommand{\li}{[\![}
\newcommand{\ri}{]\!]}
\newcommand{\loc}{\mathrm{loc}}
\newcommand{\rad}{\mathrm{rad}}
\renewcommand{\b}{\mathrm{b}}
\renewcommand{\c}{\mathrm{c}}
\renewcommand{\d}{\mathrm{d}}

\newcommand{\dist}{\mathrm{dist}}

\newcommand{\M}{\mathrm{max}}
\newcommand{\w}{{\textsl w}}
\renewcommand{\le}{\leqslant}
\renewcommand{\ge}{\geqslant}
\renewcommand{\Re}{\mathrm{Re}}
\renewcommand{\Im}{\mathrm{Im}}
\newcommand{\bs}{\boldsymbol}
\newcommand{\vi}{\bs{\mathrm{i}}}
\newcommand{\p}{\prime}
\newcommand{\vect}{\overrightarrow}
\DeclareMathOperator{\supp}{supp}

\DeclareMathOperator{\sign}{sign}

\numberwithin{equation}{section}

\newtheorem{thm}{Theorem}[section]
\newtheorem{prop}[thm]{Proposition}
\newtheorem{cor}[thm]{Corollary}
\newtheorem{lem}[thm]{Lemma}

\theoremstyle{definition}
\newtheorem{rmk}[thm]{Remark}
\newtheorem{defi}[thm]{Definition}

\newtheorem*{assumexi*}{Existence assumption}
\newtheorem*{assumuni*}{Uniqueness assumption}

\newenvironment{proof*}{\noindent{\bf Proof.}}{\qed}
\newenvironment{vproof}[1]{\noindent{\bf Proof #1}}{\qed}

\title{\huge \sc Localizing Estimates of the Support of Solutions of some Nonlinear Schrödinger Equations -- The Stationary Case}
\author{\sc Pascal Bégout\footnote{The research of Pascal Bégout was partially supported by Grants HPRN--CT--2002-00274 of the European program {\it Nonlinear partial differential equations describing front propagation and other singular phenomena}} and Jes\'us Ildefonso
D{\'{\i}}az\footnote{The research of J.I.~D\'{\i}az was partially supported by the project ref. MTM200806208 of the DGISPI (Spain) and the Research Group MOMAT (Ref. 910480) supported by UCM. He has received also support from the ITN \textit{FIRST} of the Seventh Framework Program of the European Community's (grant agreement number 238702)
}}
\date{}

\begin{document}

\maketitle

\begin{gather*}
\begin{array}{cc}
    ^*\mbox{Institut de Mathématiques de Toulouse }	& \;^\dagger\mbox{Departamento de Matem\'atica Aplicada} \\
                                   \mbox{Université de Toulouse }	& \mbox{ Facultad de Ciencias Matem\'aticas} \\
                                 \mbox{Manufacture des Tabacs }	& \mbox{ Universidad Complutense de Madrid} \\
                                        \mbox{21, Allée de Brienne }	& \mbox{ Plaza de Ciencias, 3} \\
                \mbox{31000 Toulouse Cedex, FRANCE }	& \mbox{ 28040 Madrid, SPAIN}
\bigskip \\
\mbox{
{\footnotesize$^*$e-mail\:: }\htmladdnormallink{{\footnotesize\udl{\tt{Pascal.Begout@math.cnrs.fr}}}}{mailto:Pascal.Begout@math.cnrs.fr}}
&
\mbox{
{\footnotesize$^\dagger$e-mail\:: }\htmladdnormallink{{\footnotesize\udl{\tt{ildefonso\_diaz@mat.ucm.es}}}}
{mailto:ildefonso_diaz@mat.ucm.es}
}
\end{array}
\end{gather*}

\begin{abstract}
The main goal of this paper is to study the nature of the support of the solution of suitable nonlinear Schrödinger equations, mainly the compactness of the support and its spatial localization. This question touches the very foundations underlying the derivation of the Schrödinger equation, since it is well-known a solution of a linear Schrödinger equation perturbed by a regular potential never vanishes on a set of positive measure. A fact, which reflects the impossibility of locating the particle. Here we shall prove that if the perturbation involves suitable singular nonlinear terms then the support of the solution is a compact set, and so any estimate on its spatial localization implies very rich information on places not accessible by the particle. Our results are obtained by the application of certain energy methods which connect the compactness of the support with the local vanishing of a suitable ``energy function'' which satisfies a nonlinear differential inequality with an exponent less than one. The results improve and extend a previous short presentation by the authors published in 2006.
\end{abstract}

{\let\thefootnote\relax}
{\let\thefootnote\relax}

\tableofcontents

\baselineskip .65cm

\section{Introduction}
\label{intro}

This paper deals with the study of the following stationary nonlinear Schrödinger equation (SNLS) with a complex singular potential
\begin{gather}
 \label{snls}
  -\vi\Delta\u+\va|\u|^{-(1-m)}\u+\bs{bu}=\F(x), \mbox{ in } \Omega.
\end{gather}
Here, $\Omega \subseteq \mathbb{R}^{N}$ is an open subset, $0<m<1,$ and $(\va,\vb)\in\C^2.$ The interest of the consideration of this stationary problem is motivated not only in order to study the asymptotic states, when $t\tends\infty,$ of the solutions of the associated evolution problem but also by the study of the so called standing waves of the evolution problem~\eqref{nls} below, with $\vb\in\vi\R$ in \eqref{snls}. Indeed, choosing arbitrarily $\bs b\in\vi\R$ in~\eqref{snls} and setting for any $(t,x)\in\R\times\Omega,$ $\bs\vphi(t,x)=\u(x)e^{\bs b t},$ if $\u$ is a solution
to~\eqref{snls} then $\bs\vphi$ is a solution to
\begin{gather}
 \label{nls}
  \left\{
   \begin{split}
    \vi\frac{\partial\bs\vphi}{\partial t}+\Delta\bs\vphi+\vi\va|\bs\vphi|^{-(1-m)}\bs\vphi=\vi\F(x)e^{\bs b t},	& \mbox{ in } \R\times\Omega,			\\
    \bs\vphi_{|\partial\Omega}=\0, 														& \mbox{ on } \R\times\partial\Omega,	\\
    \bs\vphi(0)=\u, 																	& \mbox{ in } \Omega.
   \end{split}
  \right.
\end{gather}
The main goal of this paper is to study the nature of the support of the solution of~\eqref{snls}: mainly its compactness and localization. Let us mention that, in our opinion, this question touches the very foundations of the derivation of the Schrödinger equation. Indeed, one of the main modifications introduced by Quantum Mechanics, with respect Classical Mechanics, is the impossibility to localize the state (position and velocity) of a particle. The solution $\u(t,x)$ is related to the probability of finding the position and momentum of particle (see, e.g. the presentation made in the text book by Strauss~\cite{MR1159712}. It is well-known that in most of the different versions of the Schrödinger equations the corresponding solution never vanishes on a subset positive measure of the domain, which reflects the impossibility of localizing the particle as mentioned above. This is the case, for instance, in case of the linear Schrödinger equation and also for some nonlinear versions where the linear equation is perturbed by a nonlinear regular potential (see, for instance, the monographs of Sulem and Sulem~\cite{MR2000f:35139} and Cazenave~\cite{MR2002047}).
\\
The main goal of this work is to show that if the linear Schrödinger equation is perturbed with suitable singular nonlinear potentials, then the support of the solution becomes a compact set and so any estimate on its spatial localization implies very rich information on places which can not
be occupied by the particle.
\\
We point out that complex potentials with certain types of singularities arise in many different situations (see, for instance, in Brezis and
Kato~\cite{MR80i:35135}, LeMesurier~\cite{MR2000k:35256} and Liskevitch and Stollmann~\cite{MR1828819}, and the references therein). We also refer the reader to the survey Belmonte-Beitia~\cite{B} in which the author supplying many references to this type of equation and many other contexts such as: semiconductors, nonlinear optics, Bose-Einstein condensation, plasma physics, molecular dynamics. Special mention is paid in this paper to the so-called Gross-Pitaevskii (corresponding to $\vb\neq\0).$
\\
In this paper, we improve some of our previous results, outlined briefly in Bégout and D\'{\i}az~\cite{MR2214595}. Moreover, we include here new estimates and generalizations. We are aware of very few other results in the literature dealing with the support of solutions of nonlinear Schrödinger equations. For instance, Rosenau and Schochet~\cite{MR2133462} propose a (one-dimensional) quasilinear Schrödinger equation in order to get solutions with compact support for each $t$ fixed. That equation and the techniques used in that paper are very different from the ones in the present work. Analogously, in a paper dated from 2008 (\cite{PhysRevLett.101.264101}), Kashdan and Rosenau consider the question of the existence (with some numerical experiments) of \ some special solutions: an one-dimensional travelling wave solution of soliton type
$\u(t,x)=A(x-\lambda t)\exp\big(\vi(\ell(x-\lambda t)+\omega t)\big)$, for the special case of $\va=\vi\gamma$ (in problem~\eqref{nls}) and $m\in (0,1).$ They also consider the two-dimensional case (now with changing propagation directions). A nonlinear term (of cubic type) is added in their equation. Those interesting results are independent of our study which also applies in the presence of some additional nonlinear terms as in the above mentioned reference.
\\
A more restricted point of view was taken in the paper by Carles and Gallo~\cite{MR2765425} where the authors prove finite time stabilization for a linear Schrödinger equations perturbed with a suitable singular nonlinear potential. In their setting, they also prove some kind of compactness of the support of the solution by means of a different energy method, but in their case the compactness occurs merely in time and not in the spatial coordinates. 
\\
We also point out that different propagation effects have been intensively studied in the literature, but most of them are related to singularities, spectral and other properties (see, for instance, Jensen~\cite{MR797050}). The question of the compactness of the support considered here is of very different nature. 
\\
In order to present our results, we shall start by indicating some very special cases which are consequences of more technical results stated later (see Theorem~\ref{thmstag} below).

\begin{thm}
\label{thmA}
Let $0<m<1,$ let $a\in\R\setminus\{0\}$ and let $b\in\R,$ $b>0.$ Let $\F\in\bs{L^\frac{m+1}{m}}(\R^N)$ with compact support. Then there exists a unique weak solution $\u\in\bs{H^1}(\R^N)\cap\bs{L^{m+1}}(\R^N)$ $($see Definition~$\ref{defsols}$ below$)$ of the problem
\begin{gather*}
-\vi\Delta\u+a|\u|^{-(1-m)}\u+\vi b\u=\F(x), \mbox{ in } \R^N.
\end{gather*}
In addition, $\u$ is compactly supported.
\end{thm}

\begin{thm}
\label{thmB}
Let $\Omega\subseteq\R^N$ be a nonempty open subset, let $0<m<1,$ let $a\in\R\setminus\{0\}$ and let $b\in\R,$ $b>0.$ Let
$\F\in\bs{L^\frac{m+1}{m}}(\Omega)$ with compact support. Assume that $\F$ is small enough in $\bs{L^\frac{m+1}{m}}(\Omega).$ Then there exists a unique weak solution $\u\in\bs{H^1_0}(\Omega)\cap\bs{L^{m+1}}(\Omega)$ $($see Definition~$\ref{defsols}$ below$)$ of the problem
\begin{gather*}
 \left\{
  \begin{split}
   -\vi\Delta\u+a|\u|^{-(1-m)}\u+\vi b\u=\F(x),	& \mbox{ in } \Omega,		\\
   \u_{|\partial\Omega}=\0, 					& \mbox{ on } \partial\Omega.
  \end{split}
 \right.
\end{gather*}
In addition, $\u$ is compactly supported in $\Omega.$
\end{thm}

\noindent
We emphasize that no sign assumption has been made on $a$ in the precedent statements. Much more general versions of our results are presented in the next section where we also include a detailed explanation of the notations used in this paper.

\section{Notations and general versions of the main results}
\label{nota}

Before stating our main results we shall indicate here some of the notations used throughout. Bold symbols are used for complex mathematics objets. For a real number $r,$ $r_+=\max\{0,r\}$ is the positive part of $r.$ We write $\vi^2=-1.$ We denote by $\ovl{\z}$ the conjugate of the complex number $\z,$ by $\Re(\z)$ its real part and by $\Im(\z)$ its imaginary part. For $1\le p\le\infty,$ $p^\p$ is the conjugate of $p$ defined by $\frac{1}{p}+\frac{1}{p^\p}=1.$ Let $j,k\in\Z$ with $j<k.$ We then write $[\![j,k]\!]=[j,k]\cap\Z.$ We denote by $\partial\Omega$ the boundary of a nonempty subset $\Omega\subseteq\R^N,$ $\ovl\Omega$ its closure, $\Omega^\c=\R^N\setminus\Omega$ its complement and
$\omega\Subset\Omega$ means that $\omega\subset\Omega$ and that $\omega$ is a compact subset of $\R^N.$ For an open subset
$\Omega\subseteq\R^N,$ the usual Lebesgue and Sobolev spaces are respectively denoted by $\bs{L^p}(\Omega)=\bs{L^p}(\Omega;\C)$ and
$\bs{W^{m,p}}(\Omega)=\bs{W^{m,p}}(\Omega;\C)$ $(1\le p\le\infty$ and $m\in\N),$ $\bs{H^m}(\Omega)=\bs{W^{m,2}}(\Omega;\C),$
$\bs{H^m_0}(\Omega)=\bs{W^{m,2}_0}(\Omega;\C)$ is the closure of $\bs{\Dr}(\Omega)=\bs{\Dr}(\Omega;\C)$ under the $\bs{H^m}$-norm, and
$\bs{H^{-m}}(\Omega)$ is its topological dual. $\bs{H^1_\c}(\Omega)=\left\{\u\in\bs{H^1}(\Omega);\supp\u\Subset\Omega\right\}.$
$\bs C(\Omega)=\bs{C^0}(\Omega)=\bs C(\Omega;\C)=\bs{C^0}(\Omega;\C)$ is the space of continuous functions from $\Omega$ to $\C.$ For
$k\in\N,$ $\bs{C^k}(\Omega)=\bs{C^k}(\Omega;\C)$ is the space of functions lying in $\bs C(\Omega;\C)$ and having all derivatives of order lesser or equal than $k$ belonging to $\bs C(\Omega;\C).$ For $0<\alpha\le1$ and $k\in\N_0\stackrel{\mathrm{def}}{=}\N\cup\{0\},$
$\bs{C^{k,\alpha}_\loc}(\Omega)=\bs{C^{k,\alpha}_\loc}(\Omega;\C)
=\left\{\u\in\bs{C^k}(\Omega;\C);\forall\omega\Subset\Omega,\;\vsum_{|\beta|=k}H^\alpha_\omega(D^\beta\u)<+\infty\right\},$ where
$H^\alpha_\omega(\u)=\vsup_{\left\{\substack{(x,y)\in\omega^2 \\ x\not=y}\right.}\frac{|\u(x)-\u(y)|}{|x-y|^\alpha}.$ The Laplacian in $\Omega$ is written $\Delta=\sum\limits_{j=1}^N\frac{\partial^2}{\partial x^2_j}.$ For a functional space $\bs E\subset\bs{L^1_\loc}(\Omega;\C),$ we denote by
$\bs{E_\rad}$ the space of functions $\f\in\bs E$ such that $\f$ is spherically symmetric. For a Banach space $E,$ we denote by $E^\star$ its topological dual and by $\langle\: . \; , \: . \:\rangle_{E^\star,E}\in\R$ the $E^\star-E$ duality product. In particular, for any
$\bs T\in\bs{L^{p^\p}}(\Omega)$ and $\bs{\vphi}\in\bs{L^p}(\Omega)$ with $1\le p<\infty,$
$\langle\bs T,\bs{\vphi}\rangle_{\bs{L^{p^\p}}(\Omega),\bs{L^p}(\Omega)}=\Re\vint_\Omega\bs T(x)\ovl{\bs{\vphi}(x)}\d x.$ For $x_0\in\R^N$ and $r>0,$ we denote by $B(x_0,r)=\{x\in\R^N;\; |x-x_0|<r\}$ the open ball of $\R^N$ of center $x_0$ and radius $r,$ by
$\S(x_0,r)=\{x\in\R^N;\; |x-x_0|=r\}$ its boundary and by $\ovl B(x_0,r)=B(x_0,r)\cup\S(x_0,r)$ its closure. We also use the notation
$B_\Omega(x_0,r)=\Omega\cap B(x_0,r).$ As usual, we denote by $C$ auxiliary positive constants, and sometimes, for positive parameters $a_1,\ldots,a_n,$ write $C(a_1,\ldots,a_n)$ to indicate that the constant $C$ continuously depends only on $a_1,\ldots,a_n$ (this convention also holds for constants which are not denoted by ``$C$'').
\\
Let us return to equation~\eqref{nls}. Note that no boundary condition is imposed since all the compact support results (which are due to
Theorem~\ref{thmstag} below) rest on the notion of local solution (Definition~\ref{defsols} below). If $\Omega\neq\R^N,$ boundary conditions are necessary for establishing existence and uniqueness of global solutions of \eqref{snls}. For the purpose of clarity, we shall consider the Dirichlet case,
\begin{gather}
 \label{dbc}
  \u_{|\partial\Omega}=\0, \mbox{ on }\partial\Omega,
\end{gather}
rather than Neumann boundary condition, mixed boundary condition or another one. The choice of the boundary condition is motivated by the integration by parts relation $\langle\Delta u,v\rangle=-\langle\nabla u,\nabla v\rangle.$

\noindent
Compactness, existence and uniqueness results will follow from assumptions on $(\va,\vb)\in\C^2$ stated below. Define the following subsets
\begin{gather*}
	\begin{cases}
		\A=\C\setminus\big\{\z\in\C; \Re(\z)=0 \text{ and } \Im(\z)\le0\big\},	\medskip \\
		\B=\A\cup\big\{\0\big\}.
	\end{cases}
\end{gather*}
\begin{assumexi*}
Let $(\va,\vb)\in\C^2$ satisfy
\begin{gather}
\label{ab}
(\va,\vb)\in\A\times\B \quad \text{ and } \quad
	\begin{cases}
		\Re(\va)\Re(\vb)\ge0,											\medskip \\
		\text{ or }														\medskip \\
		\Re(\va)\Re(\vb)<0 \; \text{ and } \; \Im(\vb)>\dfrac{\Re(\vb)}{\Re(\va)}\Im(\va).
	\end{cases}
\end{gather}
\end{assumexi*}

\begin{assumuni*}
Let $(\va,\vb)\in\C^2$ satisfy
\begin{gather}
\label{abuni}
\Im(\va)\ge0 \quad \text{ and } \quad
	\begin{cases}
		\va\neq\0 \; \text{ and } \; \Re(\va\ovl\vb)\ge0,	\medskip \\
		\text{ or }								\medskip \\
		\va=\0 \; \text{ and } \; \vb\in\B.
	\end{cases}
\end{gather}
\end{assumuni*}
\noindent
For a geometric explanation of these hypotheses, see Section~\ref{picture}. For $(\va,\vb)\in\C^2$ satisfying~\eqref{ab}, it will be convenient to introduce the following constants. Let $\delta>0$ be an arbitrarily chosen parameter. 
\begin{gather}
\label{A}
A(\delta)=\frac{|\Re(\va)|+|\Im(\va)|+\delta}{|\Re(\va)|}, \text{ if } \Re(\va)\neq0,
\end{gather}
\begin{gather}
\label{B}
B=\dfrac{|\Re(\vb)|+|\Im(\vb)|}{|\Re(\vb)|}, \text{ if } \Re(\vb)\neq0,
\end{gather}
\begin{gather}
\label{L}
L=
	\begin{cases}
		\delta,							&	\text{if }	\Im(\va)<0		\; \text{ and } \; 	\Re(\va)\Re(\vb)\ge0,			\medskip \\
		|\Re(\va)|,							&	\text{if }	\Im(\va)=0,	\; \Im(\vb)\ge0 \;	\text{ and } \; \Re(\va)\Re(\vb)\ge0,	\medskip \\
		\Im(\va)							&	\text{if }	\Im(\va)>0		\; \text{ and } \;	\Im(\vb)\ge0,					\medskip \\
		\Im(\va)-\dfrac{\Re(\va)}{\Re(\vb)}\Im(\vb),	&	\text{otherwise},
	\end{cases}
\end{gather}
\begin{gather}
\label{M}
M=
    \begin{cases}
	\max\big\{A(\delta),B\big\},		&	\text{if } 	\Im(\va)<0,	\; \Im(\vb)<0 \;	\text{ and } \;	\Re(\va)\Re(\vb)\ge0,		\medskip \\
	A(\delta),					&	\text{if }	\Im(\va)<0,	\; \Im(\vb)\ge0 \;	\text{ and } \;	\Re(\va)\Re(\vb)\ge0,		\medskip \\
	2						&	\text{if }	\Im(\va)\ge0,	\; \Im(\vb)\ge0 \;	\text{ and } \;	\big(\Im(\va)>0 \; \text{ or }
																			\; \Re(\va)\Re(\vb)\ge0\big),	\medskip \\
	B						&	\text{if }	\big(\Im(\va)\ge0	\; \text{ and } \;	\Im(\vb)<0\big)	\; \text{ or }\; \Re(\va)\Re(\vb)<0.
    \end{cases}
\end{gather}
Under hypothesis~\eqref{ab}, one easily checks that $A(\delta),$ $B,$ $L$ and $M$ are well defined and positive. The parameter $\delta$ may seem very mysterious but, actually, it is not. In order to obtain the crucial estimate \eqref{proofthmstag1}, we apply Lemma~\ref{lemAB} to
\eqref{proofthmstag2} and \eqref{proofthmstag3}. The hard case $\Im(\va)<0$ can be treated in the following way. If $\Re(\va)\Re(\vb)>0$ then we add the assumption $\Im(\vb)>\frac{\Re(\vb)}{\Re(\va)}\Im(\va).$ But when $\Re(\va)\Re(\vb)\le0,$ if we do not want make an additional assumption on $\va$ and $\vb,$ we have to introduce a positive parameter $\delta$ in order to obtain a positive coefficient $L=L(\delta)$ in front of
$\|\u\|_{\bs{L^{m+1}}(B(x_0,\rho))}^{m+1}$ (played by $C_2$ in Lemma~\ref{lemAB}). If we do not introduce this parameter (that is, if we choose
$\delta=0)$ then we get $L=0$ in \eqref{proofthmstag1} and we loose the effect of the nonlinearity (see Cases~5 and 6 in the proof of
Lemma~\ref{lemAB}).
\medskip \\
Numerical computations of stationary solutions are done in Bégout and Torri~\cite{BegTor}, while the evolution case and self-similar solutions are studied in Bégout and D\'iaz~\cite{Beg-Di2,Beg-Di3}, respectively. In this paper, we prove the results stated in Bégout and D\'iaz~\cite{MR2214595} and add some generalizations. This paper is concerned with the propagation of the support of $\F$ to the solution $\u,$ and all these results are a consequence of the following theorem.

\begin{thm}
\label{thmstag}
Let $\Omega\subseteq\R^N$ be a nonempty open subset, let $0<m<1,$ let $(\va,\vb)\in\C^2$ satisfying~\eqref{ab}, let $L>0$ be given by~\eqref{L} and let $M>0$ be given by~\eqref{M}. There exists $C=C(N,m)>0$ satisfying the following property. Let $\F\in\bs{L^1_\loc}(\Omega),$ let
$\u\in\bs{H^1_\loc}(\Omega)$ be any local weak solution of \eqref{snls} $($see Definition~$\ref{defsols}$ below$),$ let $x_0\in\Omega$ and let
$\rho_0>0.$ If $\rho_0>\dist(x_0,\partial\Omega)$ then assume further that $\u\in\bs{H^1_0}(\Omega).$ If $\F_{|B_\Omega(x_0,\rho_0)}\equiv\0$ then
$\u_{|B_\Omega(x_0,\rho_\M)}\equiv\0,$ where
\begin{multline}
 \label{thmstag1}
  \rho_\M^\nu=\left(\rho_0^\nu-CM^2\max\left\{1,\frac{1}{L^2}\right\}\max\left\{\rho_0^{\nu-1},1\right\}\right. \\
  \left.\times\vmin_{\tau\in\left(\frac{m+1}{2},1\right]}\left\{\frac{E(\rho_0)^{\gamma(\tau)}
  \max\{b(\rho_0)^{\mu(\tau)},b(\rho_0)^{\eta(\tau)}\}}{2\tau-(1+m)}\right\}\right)_+,
\end{multline}
and where for any $\tau\in\left(\frac{m+1}{2},1\right],$
\begin{gather*}
\begin{array}{lll}
     E(\rho_0)=\|\nabla\u\|_{\bs{L^2}(B_\Omega(x_0,\rho_0))}^2, &
b(\rho_0)=\|\u\|_{\bs{L^{m+1}}(B_\Omega(x_0,\rho_0))}^{m+1}, & 
                                                                                                      \gamma(\tau)=\frac{2\tau-(1+m)}{k}\in(0,1), \medskip \\
\mu(\tau)=\frac{2(1-\tau)}{k}, & \eta(\tau)=\frac{1-m}{1+m}-\gamma(\tau)>0, 	& k=2(1+m)+N(1-m),  \medskip \\
\nu=\frac{k}{m+1}>2.              &										&
\end{array}
\end{gather*}
\end{thm}

\begin{rmk}
\label{rmkthmstag}
If the solution is too ``large'', it may happen that $\rho_\M=0$ and so the above result is not consistent. A sufficient condition to observe a localizing effect is that the solution is small enough, in a suitable sense. We give two results in this direction. The first one (Theorem~\ref{thmstacom}) pertains to the size of the solution, while the second one is concerned with the size of the external source $\F$ (Theorem~\ref{thmstaO}), which seems to be more natural. In addition, Theorem~\ref{thmstaO} says where the support of the solutions is localized with respect to the support of the external source $\F.$
\end{rmk}

\noindent
Now, we state the precise notion of solution.

\begin{defi}
\label{defsols}
Let $\Omega\subseteq\R^N$ be an open subset, let $(\va,\vb)\in\bs{\C^2},$ let $0<m<1$ and let $\F\in\bs{L^1_\loc}(\Omega).$ We say that $\u$ is a {\it local weak solution} of \eqref{snls} if $\u\in\bs{H^1_\loc}(\Omega)$ and if $\u$ is a solution of \eqref{snls} in $\bs{\Dr^\p}(\Omega),$ that is
\begin{gather}
\label{defsols1}
\langle-\vi\Delta\u+\va|\u|^{-(1-m)}\u+\bs{bu},\bs\vphi\rangle_{\bs{\Dr^\p}(\Omega),\bs{\Dr}(\Omega)}
=\langle\F,\bs\vphi\rangle_{\bs{\Dr^\p}(\Omega),\bs{\Dr}(\Omega)},
\end{gather}
for any $\bs\vphi\in\bs{\Dr}(\Omega).$ \\
We say that $\u$ is a {\it global weak solution} of \eqref{snls} and \eqref{dbc} if $\u$ is a local weak solution of \eqref{snls} and if furthermore
$\u\in\bs{H^1_0}(\Omega)\cap\bs{L^{m+1}}(\Omega).$ \\
Let $\z\in\C\setminus\{\0\}.$ Since $\left||\z|^{-(1-m)}\z\right|=|\z|^m,$ it is understood that $\left||\z|^{-(1-m)}\z\right|=0$ when $\z=\0.$
\end{defi}

\begin{rmk}
\label{rmkdefsols}
Here are some comments about Definition~\ref{defsols}.
\begin{enumerate}
 \item
  \label{rmkdefsols1}
   For a global weak solution $\u$ of~\eqref{snls} and \eqref{dbc}, the boundary condition $\u_{|\partial\Omega}=0$ is included in the assumption
   $\u\in\bs{H^1_0}(\Omega).$ On the contrary, the notion of local weak solution does not consider any boundary condition.
 \item
  \label{rmkdefsols2}
   When $\u$ is a local weak solution of~\eqref{snls}, we have $\nabla\u\in\bs{L^2_\loc}(\Omega),$
   $\va|\u|^{-(1-m)}\u\in\bs{L^\frac{m+1}{m}_\loc}(\Omega)$ and $\bs{bu}\in\bs{L^2_\loc}(\Omega).$ Then $\Delta\u\in\bs{L^1_\loc}(\Omega)$ and
   equation~\eqref{snls} makes sense in $\bs{L^1_\loc}(\Omega).$ Furthermore, $\bs{L^\frac{m+1}{m}_\loc}(\Omega)\subset\bs{L^2_\loc}(\Omega)$
   and $\bs\Dr(\Omega)$ is dense in $\bs{H^1_\c}(\Omega).$ It follows from Sobolev's embedding that if $\u$ is a local weak solution of~\eqref{snls}
   then
    \begin{multline}
     \label{defsols2}
      \Re\vint_\Omega\vi\nabla\u(x).\ovl{\nabla\bs\vphi(x)}\d x+\Re\vint_\Omega\left(\va|\u(x)|^{-(1-m)}\u(x)+\bs{bu}(x)\right)\ovl{\bs\vphi(x)}\d x \\
      =\Re\vint_\Omega\F(x)\ovl{\bs\vphi(x)}\d x,
    \end{multline}
   for any $\bs \vphi\in\bs{H^1_\c}(\Omega)$ with either $\supp\bs\vphi\cap\supp\F=\emptyset$ or $\F\in\bs{L^\frac{p}{p-1}_\loc}(\Omega),$ for some
   $1\le p\le\infty$ if $N=1,$ $1\le p<\infty$ if $N=2$ or $1\le p\le\frac{2N}{N-2},$ if $N\ge3.$ For example, $p=m+1$ is always an admissible value.
 \item
  \label{rmkdefsols3}
   In the same way, by density of $\bs\Dr(\Omega)$ in $\bs{H^1_0}(\Omega)\cap\bs{L^{m+1}}(\Omega)\cap\bs{L^p}(\Omega),$ for any $1\le p<\infty,$
   and in $\bs{H^1_0}(\Omega)\cap\bs{L^{m+1}}(\Omega),$ if $\u$ is a global weak solution of \eqref{snls} and \eqref{dbc} then \eqref{defsols2} holds
   for any $\bs\vphi\in\bs{H^1_0}(\Omega)\cap\bs{L^{m+1}}(\Omega)$ with either $\supp\bs\vphi\cap\supp\F=\emptyset$ or
   $\bs\vphi\in\bs{L^p}(\Omega)$ and $\F\in\bs{L^\frac{p}{p-1}}(\Omega),$ for some $1\le p<\infty.$ In particular, if $p$ is as in~\ref{rmkdefsols2}.~of
   this remark with additionally $p\ge m+1,$ then in view of $\bs{H^1_0}(\Omega)\cap\bs{L^{m+1}}(\Omega)\inj\bs{L^p}(\Omega),$
   equation~\eqref{snls} makes sense in $\bs{H^{-1}}(\Omega)+\bs{L^\frac{m+1}{m}}(\Omega)$ and \eqref{defsols2} holds for any
   $\bs\vphi\in\bs{H^1_0}(\Omega)\cap\bs{L^{m+1}}(\Omega).$
\end{enumerate}
\end{rmk}

\section{Spatial localization property}
\label{staresults}

\begin{thm}
\label{thmF}
Let $\Omega\subseteq\R^N$ be a nonempty open subset, let $0<m<1$ and let $(\va,\vb)\in\C^2$ satisfying \eqref{ab}. Let
$\F\in\bs{L^\frac{m+1}{m}}(\Omega),$ let $\u\in\bs{H^1_\loc}(\Omega)$ be any local weak solution of \eqref{snls} $($Definition~$\ref{defsols}),$ let $x_0\in\Omega$ and let $\rho_1>0.$ If $\rho_1>\dist(x_0,\partial\Omega)$ then assume further that $\u\in\bs{H^1_0}(\Omega).$ Then there exist
$E_\star>0$ and $\eps_\star>0$ satisfying the following property. Let $\rho_0\in(0,\rho_1).$ If
$\|\nabla\u\|_{\bs{L^2}(B_\Omega(x_0,\rho_1))}^2<E_\star$ and
\begin{gather}
\label{thmF1}
\forall\rho\in(0,\rho_1), \;
\|\F\|_{\bs{L^\frac{m+1}{m}}(B_\Omega(x_0,\rho))}^\frac{m+1}{m}\le\eps_\star\big((\rho-\rho_0)_+\big)^p,
\end{gather}
where $p=\frac{2(1+m)+N(1-m)}{1-m}>N+2,$ then $\u_{|B_\Omega(x_0,\rho_0)}\equiv\0.$ In other words, with the notation of
Theorem~$\ref{thmstag},$ $\rho_\M=\rho_0.$
\end{thm}

\begin{rmk}
\label{rmkthmF}
We may estimate $E_\star$ and $\eps_\star$ as
\begin{gather*}
E_\star=E_\star\left(\|\u\|_{\bs{L^{m+1}}(B(x_0,\rho_1))}^{-1},\rho_1,\frac{\rho_0}{\rho_1},\frac{L}{M},N,m\right), \\
\eps_\star=\eps_\star\left(\|\u\|_{\bs{L^{m+1}}(B(x_0,\rho_1))}^{-1},\frac{\rho_0}{\rho_1},\frac{L}{M},N,m\right),
\end{gather*}
where $L>0$ and $M>0$ are given by~\eqref{A} and~\eqref{M}, respectively. The dependence on $\frac{1}{\delta}$ means that for any value $\delta$ small enough, $E_\star$ and $\eps_\star$ are bounded from below. \\
Note that $p=\frac{1}{\gamma(1)},$ where $\gamma$ is the function defined in Theorem~\ref{thmstag}.
\end{rmk}

\begin{thm}
\label{thmstacom}
Let $\Omega\subseteq\R^N$ be a nonempty open subset, let $0<m<1,$ let $(\va,\vb)\in\C^2$ satisfying~\eqref{ab}, let $L>0$ be given by~\eqref{L} and let $M>0$ be given by~\eqref{M}. There exists $C=C(N,m)>0$ satisfying the following property. Let $\bs{F\in L^1_\loc}(\Omega),$ let $\u\in\bs{H^1_\loc}(\Omega)$ be any local weak solution of \eqref{snls} $($Definition~$\ref{defsols}),$ let $x_0\in\Omega$ and let $\rho_0>0.$ If
$2\rho_0>\dist(x_0,\partial\Omega)$ then assume further that $\u\in\bs{H^1_0}(\Omega).$ Finally, suppose
$\F_{|B_\Omega(x_0,2\rho_0)}\equiv\0,$ $\|\u\|_{\bs{L^{m+1}}(B_\Omega(x_0,2\rho_0))}\le1$ and one of the two estimates \eqref{thmstacom1} or \eqref{thmstacom2} below is satisfied.
\begin{gather}
 \label{thmstacom1}
  \|\nabla\u\|_{\bs{L^2}(B_\Omega(x_0,2\rho_0))}^\frac{2(1-m)}{k}\le
  C(2^\nu-1)(1-m)M^{-2}\min\left\{1,L^2\right\}\min\left\{\dfrac{1}{2},\rho_0\right\}^{\nu-1}\rho_0, \medskip \\
 \begin{cases}
  \label{thmstacom2}
   \|\nabla\u\|_{\bs{L^2}(B_\Omega(x_0,2\rho_0))}\le1, \medskip \\
   \|\u\|_{\bs{L^{m+1}}(B_\Omega(x_0,2\rho_0))}^\frac{2s(m+1)}{k}\le
   C(2^\nu-1)(1-m-2s)M^{-2}\min\left\{1,L^2\right\}\min\left\{\dfrac{1}{2},\rho_0\right\}^{\nu-1}\rho_0,
 \end{cases}
\end{gather}
for some $s\in\left(0,\frac{1-m}{2}\right),$ where the constants $k>\nu>2$ are given in Theorem~$\ref{thmstag}.$ Then
$\u_{|B_\Omega(x_0,\rho_0)}\equiv\0.$
\end{thm}

\begin{rmk}
\label{rmkthmstacom}
Note that in estimate~\eqref{thmstacom1}, $\frac{2(1-m)}{k}=\frac2p,$ where $p>N+2$ is given in Theorem~\ref{thmF}.
\end{rmk}

\begin{thm}
\label{thmstaO}
Let $\Omega\subseteq\R^N$ be a nonempty open subset, let $0<m<1,$ let $(\va,\vb)\in\C^2$ satisfying \eqref{ab}, let $L>0$ be given by~\eqref{L} and let $M>0$ be given by~\eqref{M}. Then for any $\eps>0,$ there exists $\delta_0=\delta_0(\eps,N,m,L,M)>0$ satisfying the following property. Let $\F\in\bs{L^\frac{m+1}{m}}(\Omega)$ and let $\u\in\bs{H^1_0}(\Omega)\cap\bs{L^{m+1}}(\Omega)$ be any global weak solution of \eqref{snls} and \eqref{dbc}. If $\supp\F$ is a compact set and if $\|\F\|_{\bs{L^\frac{m+1}{m}}(\Omega)}\le\delta_0$ then
$\supp\u\subset\ovl\Omega\cap\vO(\eps),$ where $\vO(\eps)$ is the open bounded set
\begin{gather*}
 \vO(\eps)=\left\{x\in\R^N;\; \exists y\in\supp\F \;\; such \;\; that \;\; |x-y|<\eps \right\}.
\end{gather*}
In particular, if $\eps>0$ is small enough then $\supp\u\subset\vO(\eps)\subset\Omega.$
\end{thm}

\noindent
We see that localization effect occurs under some smallness condition, either on the solution $\u$ (Theorem~\ref{thmstacom}) or on the external source $\F$ (Theorem~\ref{thmstaO}). When $\Omega=\R^N,$ the phenomenon is simpler since localization effect is always observed, without any condition of the size, neither on the solution nor on the external source, as show the following result.

\begin{thm}
\label{thmcompRN}
Let $0<m<1,$ let $(\va,\vb)\in\C^2$ satisfying \eqref{ab}, let $\F\in\bs{L^p}(\R^N),$ for some
$1\le p\le\infty,$ and let $\u\in\bs{H^1}(\R^N)\cap\bs{L^{m+1}}(\R^N)$ be any global weak solution of \eqref{snls}. If $\supp\F$ is a compact set then $\supp\u$ is also compact.
\end{thm}

\section{Existence and smoothness}
\label{es}

In this section, we give an existence result of solutions for equation~\eqref{snls} (Theorem~\ref{thmexista}), some \textit{a priori} bounds for the solutions of equation~\eqref{snls} (Theorem~\ref{bound}), which will be useful to establish our existence result, and a smoothness result for equation~\eqref{snls} (Proposition~\ref{propregedp}).

\begin{thm}
\label{thmexista}
Let $\Omega\subseteq\R^N$ be a nonempty open subset, let $0<m<1,$ let $(\va,\vb)\in\C^2$ satisfying \eqref{ab} and let
$\F\in\bs{L^\frac{m+1}{m}}(\Omega).$ Then equations \eqref{snls} and \eqref{dbc} admits at least one global weak solution
$\u\in\bs{H^1_0}(\Omega)\cap\bs{L^{m+1}}(\Omega).$ Furthermore, the following properties hold for any global weak solution $\u$ $\big($except Property~$3)\big).$
\begin{enumerate}
 \item[]
  \begin{enumerate}[$1)$]
   \item
    $\u\in\bs{W_\loc^{2,\frac{m+1}{m}}}(\Omega).$
   \item
    Let $\alpha\in(0,m].$ If $\F\in\bs{C_\loc^{0,\alpha}}(\Omega)$ then $\u\in\bs{C_\loc^{2,\alpha}}(\Omega).$
   \item
    If $\Omega=\left\{x\in\R^N;\;r<|x|<R\right\},$ for some $-\infty<r\le r_+<R\le+\infty,$ and if $\F$ is spherically
    symmetric then there exists a spherically symmetric global weak solution
    $\u\in\bs{H^1_0}(\Omega)\cap\bs{L^{m+1}}(\Omega)$ of \eqref{snls} and \eqref{dbc}. For $N=1,$ this means that if
    $\F$ is an even $($respectively, an odd$)$ function on $\Omega=(-R,-r)\cup(r,R)$ then $\u$ is also an even
    $($respectively, an odd$)$ function.
  \end{enumerate}
\end{enumerate}
\end{thm}

\begin{rmk}
\label{rmkthmess}
Assume $\F$ is spherically symmetric. Since we do not know, in general, if we have uniqueness of the solution, we are not able to show that any solution is radially symmetric. For a uniqueness result, see Theorem~\ref{thmuni} below.
\end{rmk}

\begin{rmk}
\label{rmkthmfin}
Assume $|\Omega|<\infty.$ There exists $\eps=\eps(N)>0$ such that for any $(\va,\vb)\in\C^2,$ $0<m<1$ and $\F\in\bs{L^2}(\Omega),$ if
$|\vb||\Omega|^\frac{2}{N}<\eps$ then equations \eqref{snls} and \eqref{dbc} admits at least one global weak solution $\u\in\bs{H^1_0}(\Omega).$ In addition, $\u\in\bs{H_\loc^2}(\Omega).$ Finally, Properties~2) and 3) of Theorem~\ref{thmexista} hold. For more details, see~Bégout and
Torri~\cite{BegTor}.
\end{rmk}

\begin{thm}
\label{bound}
Let $\Omega\subseteq\R^N$ be a nonempty open subset, let $0<m<1,$ let $(\va,\vb)\in\C^2$ satisfying \eqref{ab}, let $L>0$ be given by~\eqref{L}, let $M>0$ be given by \eqref{M} and let $\F\in\bs{L^\frac{m+1}{m}}(\Omega).$ Let $\u\in\bs{H^1_0}(\Omega)\cap\bs{L^{m+1}}(\Omega)$ be any global weak solution of \eqref{snls} and \eqref{dbc}. Then we have the following estimates.
\begin{gather}
\label{aprioribound1}
\|\nabla\u\|_{\bs{L^2}(\Omega)}^2+\|\u\|_{\bs{L^{m+1}}(\Omega)}^{m+1}<M_0\|\F\|_{\bs{L^\frac{m+1}{m}}(\Omega)}^\frac{m+1}{m}, \medskip \\
\label{aprioribound2}
\|\u\|_{\bs{H^1_0}(\Omega)}^2+\|\u\|_{\bs{L^{m+1}}(\Omega)}^{m+1}
<C\wt{M_0}\left(1+\|\F\|_{\bs{L^\frac{m+1}{m}}(\Omega)}^\frac{\delta(m+1)}{m}\right)\|\F\|_{\bs{L^\frac{m+1}{m}}(\Omega)}^\frac{m+1}{m},
\end{gather}
where $M_0=M\left(\frac{2M}{L}\right)^\frac{1}{m}\max\left\{1,\frac{2}{L}\right\},$ $\delta=\frac{2(1-m)}{(N+2)-m(N-2)},$
$\wt{M_0}=M_0(1+M_0^\delta)$ and $C=C(N,m).$
\end{thm}

\begin{prop}
\label{propregedp}
Let $\va\in\C,$ let $0<m<1,$ let $\bs V\in\bs{L^r_\loc}(\Omega;\C),$ for any $1<r<\infty,$ let
$\F\in\bs{L^1_\loc}(\Omega;\C)$ and, for some $\eps>0,$ let $\u\in\bs{L^{1+\eps}_\loc}(\Omega;\C)$ $\big(\u\in\bs{L^1_\loc}(\Omega;\C)$ suffices if
$\bs V\in\bs{L^\infty_\loc}(\Omega;\C)\big)$ be a solution to
\begin{gather}
 \label{F}
  -\Delta\u + \bs{Vu} + \va|\u|^{-(1-m)}\u = \F(x), \mbox{ in } \bs{\Dr^\p}(\Omega).
\end{gather}
Let $1<q<\infty$ and suppose $\u\in\bs{L_\loc^q}(\Omega).$ Then the following regularity results hold.
\begin{enumerate}
 \item[]
  \begin{enumerate}[$1)$]
   \item
    If for some $p\in[q,\infty),$ $\F\in\bs{L^p_\loc}(\Omega)$ then $\u\in\bs{W_\loc^{2,p}}(\Omega).$
   \item
    Let $\alpha\in(0,m].$ If $(\F,\bs V)\in\bs{C_\loc^{0,\alpha}}(\Omega)\times\bs{C_\loc^{0,\alpha}}(\Omega)$ then
    $\u\in\bs {C_\loc^{2,\alpha}}(\Omega).$
  \end{enumerate}
\end{enumerate}
\end{prop}

\begin{rmk}
\label{rmkpropedp1}
Since $0<m<1$ and $\u\in\bs{L_\loc^1}(\Omega),$ one has $\bs{L_\loc^\frac{1}{m}}(\Omega)\subset\bs{L_\loc^1}(\Omega)$ and so
$|\u|^{-(1-m)}\u\in\bs{L_\loc^1}(\Omega).$ In addition, from Hölder's inequality $\bs V\u\in\bs{L^1_\loc}(\Omega)$ and it follows that
$\Delta\u\in\bs{L_\loc^1}(\Omega).$ In conclusion, equation~\eqref{F} makes senses in $\bs{L_\loc^1}(\Omega).$
\end{rmk}

\begin{rmk}
\label{rmkpropedp2}
We only state a local smoothness result since we are interested by compactly supported solutions. In this case, global smoothness is immediate. Nevertheless, one may wonder what happens when a solution is not compactly supported. We use the notation of Proposition~\ref{propregedp} and assume further that $\Omega$ is bounded\footnote{Actually, assumptions on $\Omega$ we use in this remark are $\partial\Omega$ bounded and $|\Omega|<\infty.$ But these two conditions imply that $\Omega$ is bounded.} and has a $C^{1,1}$ boundary. Let the assumptions of
Proposition~\ref{propregedp} be fulfilled and let $\u\in\bs{L^q}(\Omega),$ for some $1<q<\infty,$ be a solution
to~\eqref{F} such that $\u_{|\partial\Omega}=\0$ in the sense of the trace\footnote{Let
$\bs T:\u\tends\left\{\bs{\gamma}\u,\bs{\gamma}\frac{\partial\u}{\partial\nu}\right\}$ be the trace function defined on
$\bs{\Dr}(\ovl\Omega),$ let $1<p<\infty$ and let
$\bs{X_p}(\Omega)=\big\{\u\in\bs{L^p}(\Omega);\Delta\u\in\bs{L^p}(\Omega)\big\}.$ By density of
$\bs{\Dr}(\ovl\Omega)$ in $\bs{X_p}(\Omega),$ $\bs T$ has a continuous and linear extension from
$\bs{X_p}(\Omega)$ into $\bs{W^{-\frac{1}{p},p}}(\partial\Omega)\times\bs{W^{-1-\frac{1}{p},p}}(\partial\Omega)$ (Hörmander~\cite{MR0106333}, Theorem~2 p.503; Lions and Magenes~\cite{MR0146525}, Lemma~2.2 and Theorem~2.1 p.147; Lions and Magenes~\cite{MR0146526}, Propositions~9.1, Proposition~9.2 and Theorem~9.1
p.82; Grisvard~\cite{MR775683}, p.54). Since $\u\in\bs{L^{m+1}}(\Omega),$ it follows from equation~\eqref{F} and Hölder's inequality that $\u\in\bs{X_p}(\Omega),$ for any $1<p<m+1.$ Then ``$\u_{|\partial\Omega}=\0$ in the sense of the trace'' makes sense and means that $\bs{\gamma}\u=\0.$}.
\begin{enumerate}
 \item
  \label{rmkpropedp21}
   If for some $p\in[q,\infty),$ $\F\in\bs{L^p}(\Omega)$ and $\bs V\in\bs{L^r}(\Omega),$ $\forall r\in(1,\infty),$ then
   $\u\in\bs{W^{2,p}}(\Omega)\cap\bs{W^{1,p}_0}(\Omega).$ Indeed, recalling that if for some $1<p<\infty,$ a function
   $\v\in\bs{L^p}(\Omega)$ satisfies $\Delta\v\in\bs{L^p}(\Omega)$ and $\v_{|\partial\Omega}=\0$ in the sense of the
   trace$^2$ then $\v\in\bs{W^{2,p}}(\Omega)\cap\bs{W^{1,p}_0}(\Omega)$ (Grisvard~\cite{MR775683},
   Corollary~2.5.2.2 p.131). We then apply the bootstrap method of the proof of Proposition~\ref{propregedp} to prove
   the result, where we use the embedding $\bs{L^r}(\Omega)\inj \bs{L^s}(\Omega),$ which holds for any
   $r\ge s$ (since $\Omega$ is bounded) and the global regularity result of Grisvard~\cite{MR775683}
   (Corollary~2.5.2.2 p.131) in place of a local regularity result (Cazenave~\cite{caz-sle}, Theorem~4.1.2 p.101--102).
 \item
  \label{rmkpropedp22}
   Let $\alpha\in(0,m].$ If $\Omega$ has a $C^{2,\alpha}$ boundary and $(\F,\bs V)\in\bs{C^{0,\alpha}}(\ovl\Omega)\times
   \bs{C^{0,\alpha}}(\ovl\Omega)$ then $\u\in\bs {C^{2,\alpha}}(\ovl\Omega)\cap\bs{C_0}(\Omega)
   \footnote{For $k\in\N_0$ and $0<\alpha\le1,$ $\bs{C^{k,\alpha}}(\ovl\Omega)
   =\Big\{\u\in\bs{C^k}(\ovl\Omega;\C);\vsum_{|\beta|=k}H^\alpha_\Omega(D^\beta\u)<+\infty\Big\}
   \subset\bs{W^{k,\infty}}(\Omega)$ (since $\Omega$ is bounded) and
   $\bs{C_0}(\Omega)=\big\{\u\in\bs C(\ovl\Omega);\forall x\in\partial\Omega,\:\u(x)=\0\big\}.$}.$ Indeed, it follows
   from the above remark that $\u\in\bs{W^{2,N+1}}(\Omega)\cap\bs{H^1_0}(\Omega)$ and by Sobolev's
   embedding, $\u\in\bs{C^{0,1}}(\ovl\Omega).$ Setting
    \begin{gather*}
     \f=\F(x)-\bs{Vu}-\va|\u|^{-(1-m)}\u,
    \end{gather*}
   it then follow from equation~\eqref{F} and estimate~\eqref{propregedp5} below that
   $\f\in\bs{C^{0,\alpha}}(\ovl\Omega).$ Let
   $\v\in\bs\vC\stackrel{\text{def}}{=}\bs{C^{2,\alpha}}(\ovl\Omega)\cap\bs{C_0}(\ovl\Omega)$ be a solution to
    \begin{gather}
     \label{rmkpropedp221}
      -\Delta\vw=\f,
    \end{gather}
   given by Gilbarg and Trudinger~\cite{MR1814364}, Theorem~6.14 p.107. Since $\u\in\bs{H^1_0}(\Omega)$ is also
   a solution to~\eqref{rmkpropedp221}, uniqueness for equation~\eqref{rmkpropedp221} holds in
   $\bs{H^1_0}(\Omega)$ (Lax-Milgram's Theorem) and $\bs\vC\subset\bs{H^1_0}(\Omega),$ we conclude that
   $\u=\v$ and so $\u\in\bs\vC.$
\end{enumerate}
\end{rmk}

\noindent
We end this section by giving a result for the evolution equation (in a particular case).
\begin{cor}
\label{corevo}
Let $0<m<1,$ let $(\bs\lambda,b)\in\C\times\R$ satisfying $\bs\lambda\neq\0$ and $b\ge0.$ If $\Im(\bs\lambda)=0$ then assume further
$\Re(\bs\lambda)\le0.$ Finally, let $\F\in\bs{C^{0,m}}(\R^N)$ be compactly supported. Then there exists a solution
$\u\in\bs{C^\infty}\big(\R;\bs{C_\b^{2,m}}(\R^N)\big)$ to
\begin{gather}
 \label{nls*}
  \left\{
   \begin{split}
    \vi\frac{\partial\u}{\partial t}+\Delta\u+\bs\lambda|\u|^{-(1-m)}\u=\F(x)e^{\vi b t},	& \; \mbox{ in } \R\times\R^N,	\\
    \u(0)=\bs\vphi, 													& \; \mbox{ in } \R^N.
   \end{split}
  \right.
\end{gather}
given by
\begin{gather}
\label{sw}
\forall(t,x)\in\R\times\R^N, \; \u(t,x)=\bs\vphi(x)e^{\vi b t},
\end{gather}
where $\bs\vphi\in\bs{C_\b^{2,m}}(\R^N)$ is a solution compactly supported of
\begin{gather}
\label{bs}
-\Delta\bs\vphi-\bs\lambda|\bs\vphi|^{-(1-m)}\bs\vphi+b\bs\vphi=-\F(x), \mbox{ in } \R^N,
\end{gather}
given by Theorem~$\ref{thmexista}.$ Furthermore, for any $t\in\R,$ $\supp\u(t)$ is compact.
\end{cor}

\section{Uniqueness}
\label{uni}

\begin{thm}
\label{thmdepcon}
Let $\Omega\subseteq\R^N$ be a nonempty open subset, let $0<m<1,$ let $(\va,\vb)\in\C^2\setminus\{(\0,\0)\}$ satisfying \eqref{abuni} and let $\bs{F_1},\bs{F_2}\in\bs{L^1_\loc}(\Omega)$ be such that
$\bs{F_1}-\bs{F_2}\in\bs{L^2}(\Omega).$ Let $\bs{u_1},\bs{u_2}\in\bs{H^1_0}(\Omega)\cap\bs{L^{m+1}}(\Omega)$ be two global weak solutions of
\begin{gather}
 \label{F1}
  -\vi\Delta\bs{u_1}+\va|\bs{u_1}|^{-(1-m)}\bs{u_1}+\vb\bs{u_1}=\bs{F_1}(x), \mbox{ in } \Omega, \\
 \label{F2}
  -\vi\Delta\bs{u_2}+\va|\bs{u_2}|^{-(1-m)}\bs{u_2}+\vb\bs{u_2}=\bs{F_2}(x), \mbox{ in } \Omega,
\end{gather}
respectively. We have the following estimates.
\begin{gather}
\begin{cases}
\label{thmdepcon1}
\|\bs{u_1}-\bs{u_2}\|_{\bs{L^2}(\Omega)}
		\le\dfrac{|\va|}{\Re\left(\va\ovl\vb\right)}\|\bs{F_1}-\bs{F_2}\|_{\bs{L^2}(\Omega)},
	& \text{if } \va\neq\0 \text{ and } \Re\left(\va\ovl\vb\right)>0, \medskip \\
\|\bs{u_1}-\bs{u_2}\|_{\bs{L^2}(\Omega)}\le\dfrac{1}{b_0}\|\bs{F_1}-\bs{F_2}\|_{\bs{L^2}(\Omega)},
	& \text{if } \va=\0,
\end{cases}
\end{gather}
where $b_0=|\Re(\bs{b})|,$ if $\Re(\bs{b})\neq0$ and $b_0=\Im(\bs{b}),$ if $\Re(\bs{b})=0.$ If $\va\neq\0$ and
$\Re\left(\va\ovl\vb\right)=0$ then assume further that $\bs{u_1},\bs{u_2}\in\bs{L^\infty}(\Omega).$ Then there exists a positive constant $C=C(N,m)$ such that
\begin{gather}
\label{thmdepcon2}
\|\bs{u_1}-\bs{u_2}\|_{\bs{L^2}(\Omega)}\le
C\frac{\left(\|\bs{u_1}\|_{\bs{L^\infty}(\Omega)}+\|\bs{u_2}\|_{\bs{L^\infty}(\Omega)}\right)^{1-m}}{|\va|}
\|\bs{F_1}-\bs{F_2}\|_{\bs{L^2}(\Omega)}.
\end{gather}
\end{thm}

\begin{thm}
\label{thmuni}
Let $\Omega\subseteq\R^N$ be a nonempty open subset, let $0<m<1,$ let $(\va,\vb)\in\C^2$ satisfying \eqref{abuni} and let
$\F\in\bs{L^1_\loc}(\Omega).$ Then equations \eqref{snls} and \eqref{dbc} admit at most one global weak solution
$\u\in\bs{H^1_0}(\Omega)\cap\bs{L^{m+1}}(\Omega).$
\end{thm}

\begin{cor}
\label{corthmexiuni}
Let $\Omega\subseteq\R^N$ be a nonempty open subset, let $0<m<1,$ let $(\va,\vb)\in\A\times\B$ satisfying \eqref{abuni} and let
$\F\in\bs{L^\frac{m+1}{m}}(\Omega).$ Then equations \eqref{snls} and \eqref{dbc} admit a unique global weak solution
$\u\in\bs{H^1_0}(\Omega)\cap\bs{L^{m+1}}(\Omega).$ Furthermore, this solution satisfies Properties $1)-3)$ of Theorem~$\ref{thmexista}.$
\end{cor}

\begin{cor}
\label{corthmuni}
Let $\Omega\subseteq\R^N$ be a nonempty open subset, let $0<m<1$ and let $(\va,\vb)\in\C^2$ satisfying \eqref{abuni}. Then the problem
\begin{gather*}
 \begin{cases}
  -\vi\Delta\u+\va|\u|^{-(1-m)}\u+\vb\u=\0, \mbox{ in } \Omega, \medskip \\
  \u\in\bs{H^1_0}(\Omega)\cap\bs{L^{m+1}}(\Omega),
 \end{cases}
\end{gather*}
has for unique solution $\u\equiv\0.$
\end{cor}

\begin{cor}
\label{corfin}
Let $0<m<1,$ let $(\va,\vb)\in\A\times\B$ satisfying \eqref{abuni} and let $\F\in\bs{C^{0,m}}(\R^N)$ be compactly supported. Then there exists a unique solution $\u\in\bs{C_\b^{2,m}}(\R^N)$ of \eqref{snls} and \eqref{dbc} compactly supported. If furthermore $\F$ is spherically symmetric then
$\u$ is also spherically symmetric. For $N=1,$ this means that if $\F$ is an even $($respectively, an odd$)$ function then $\u$ is also an even
$($respectively, an odd$)$ function.
\end{cor}

\section{Pictures}
\label{picture}

In this section, we give some geometric interpretation of the values of $\va$ and $\vb.$ For convenience, we repeat the hypotheses \eqref{ab} and \eqref{abuni}. We recall that,
\begin{gather*}
	\begin{cases}
		\A=\C\setminus\big\{\z\in\C; \Re(\z)=0 \text{ and } \Im(\z)\le0\big\},	\medskip \\
		\B=\A\cup\big\{\0\big\}.
	\end{cases}
\end{gather*}
For existence of solutions to problem \eqref{snls} and \eqref{dbc}, we suppose $(\va,\vb)\in\C^2$ satisfies
\begin{gather}
\label{hypexi}
(\va,\vb)\in\A\times\B \quad \text{ and } \quad
	\begin{cases}
		\Re(\va)\Re(\vb)\ge0,											\medskip \\
		\text{ or }														\medskip \\
		\Re(\va)\Re(\vb)<0 \; \text{ and } \; \Im(\vb)>\dfrac{\Re(\vb)}{\Re(\va)}\Im(\va),
	\end{cases}
\end{gather}
while for uniqueness, we assume
\begin{gather}
\label{hypuni}
\Im(\va)\ge0 \quad \text{ and } \quad
	\begin{cases}
		\va\neq\0 \; \text{ and } \; \Re(\va\ovl\vb)\ge0,	\medskip \\
		\text{ or }								\medskip \\
		\va=\0 \; \text{ and } \; \vb\in\B.
	\end{cases}
\end{gather}
\noindent
\textbf{Existence.}
Condition~\eqref{hypexi} may easily be interpreted in this way: if $\vb\neq\0$ then one requires that $[\va,\vb]\cap\Br=\emptyset,$ where $\Br$ is the geometric representation of $\B.$ See Figures~\ref{existence1} and \ref{existence2} below. \\
\textbf{Uniqueness.}
The second condition of~\eqref{hypuni} is trivial. Indeed, $\vb$ can be chosen anywhere in the complex plane, except on the half-axis where
$\Im(\z)<0.$ Let us consider the first condition. We first choose $\va\in\C\setminus\{\0\}$ such that $\Im(\va)\ge0,$ and we choose $\vb$ with respect to $\va.$ We see $\va$ and $\vb$ as vectors of $\R^2.$ Then we write,
$
\vect a=
 \left(
  \begin{array}{c}
   \Re(\va) \\
   \Im(\va)
  \end{array}
 \right),
$
$
\vect b=
 \left(
  \begin{array}{c}
   \Re(\vb) \\
   \Im(\vb)
  \end{array}
 \right)
$
and we have
\begin{gather}
\label{abgeo}
\Re\left(\va\ovl\vb\right)=\Re(\va)\Re(\vb)+\Im(\va)\Im(\vb)=\vect a.\vect b,
\end{gather}
where $.$ denotes the scalar product between two vectors of $\R^2.$ Then the condition $\Re\left(\va\ovl\vb\right)\ge0$ is equivalent to
$\left|\angle(\vect a,\vect b)\right|\le\dfrac{\pi}{2}\rad$ (see Figure~\ref{uniqueness} below).

\begin{rmk}
\label{rmkuniexi}
Let $(\va,\vb)\in\C^2.$ Thanks to~\eqref{abgeo}, the following assertions are equivalent.
\begin{enumerate}
\item[]
\begin{enumerate}[1)]
	\item $(\va,\vb)$ satisfies \eqref{hypexi}--\eqref{hypuni} (or \eqref{ab}--\eqref{abuni}).
	\item $(\va,\vb)\in\A\times\B$ satisfies \eqref{hypuni} (or \eqref{abuni}).
	\item $\Big((\va,\vb)$ satisfies \eqref{hypuni}\Big), $\Big(\va\neq\0\Big)$ and
		$\Big(\Im(\va)=\Re(\vb)=0\implies\Im(\vb)\ge0\Big).$
\end{enumerate}
\end{enumerate}
In other words, when $\Im(\va)\neq0,$ uniqueness hypothesis~\eqref{hypuni} implies existence
hypothesis~\eqref{hypexi} (see Figure~\ref{exiuni} below).
\end{rmk}

\begin{figure}[!ht]
 \begin{minipage}[t]{7.5cm}
  \begin{center}
  \hspace{-0.6cm}
   \includegraphics[ext=pdf,keepaspectratio=true,hiresbb=true,width=8cm]{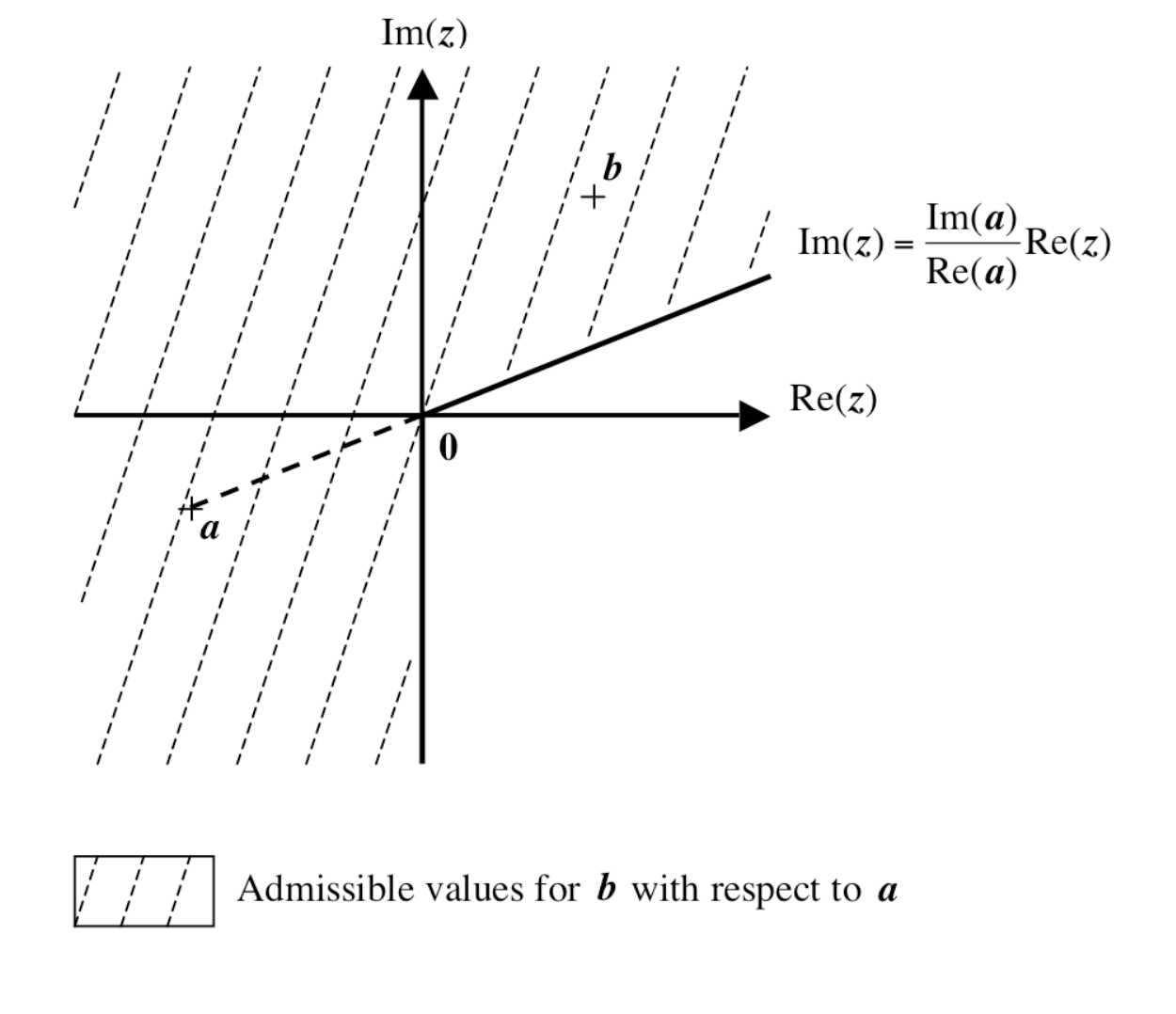}
    \caption{Existence, choice of $\vb$}
    \label{existence1}
  \end{center}
 \end{minipage} \hfill
 \begin{minipage}[t]{7.5cm}
  \begin{center}
  \hspace{-1cm}
   \includegraphics[ext=pdf,keepaspectratio=true,hiresbb=true,width=8cm]{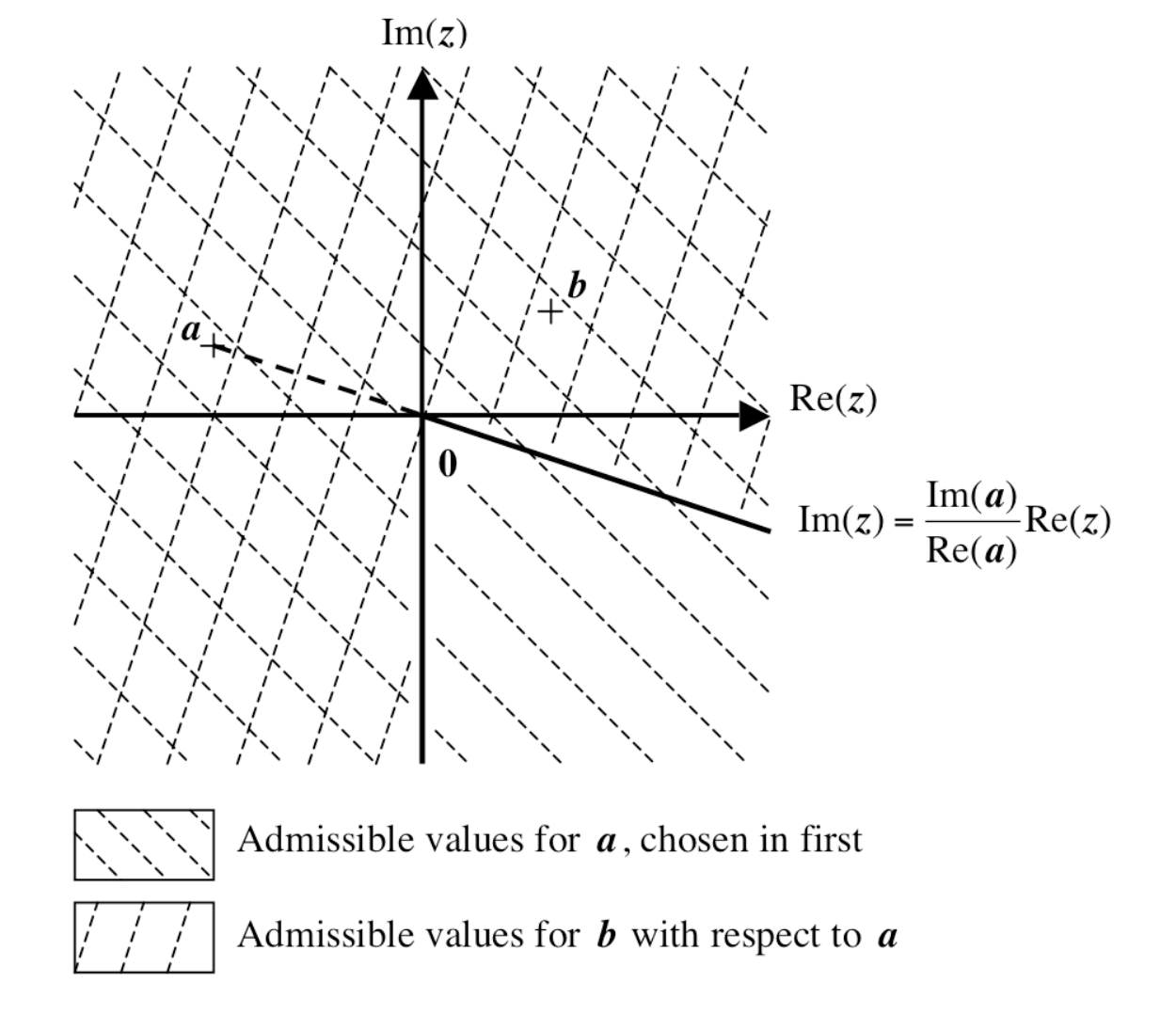}
   \caption{Existence, choice of $\va$ and $\vb$}
   \label{existence2}
  \end{center}
 \end{minipage}
 \begin{minipage}[t]{7.5cm}
  \begin{center}
  \hspace{-0.6cm}
   \includegraphics[ext=pdf,keepaspectratio=true,hiresbb=true,width=8cm]{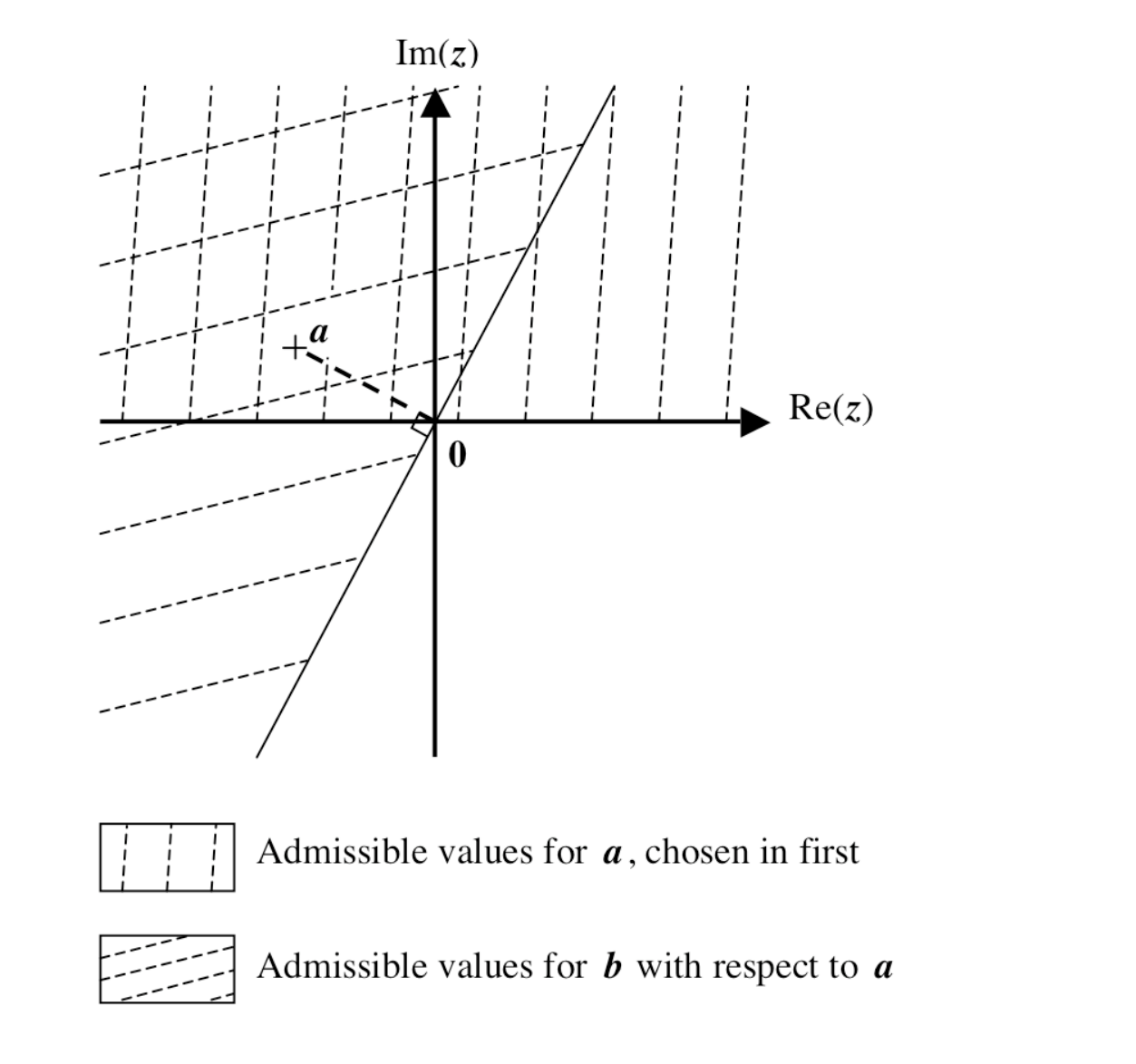}
    \caption{Uniqueness}
    \label{uniqueness}
  \end{center}
 \end{minipage} \hfill
 \begin{minipage}[t]{7.5cm}
  \begin{center}
  \hspace{-1cm}
   \includegraphics[ext=pdf,keepaspectratio=true,hiresbb=true,width=8cm]{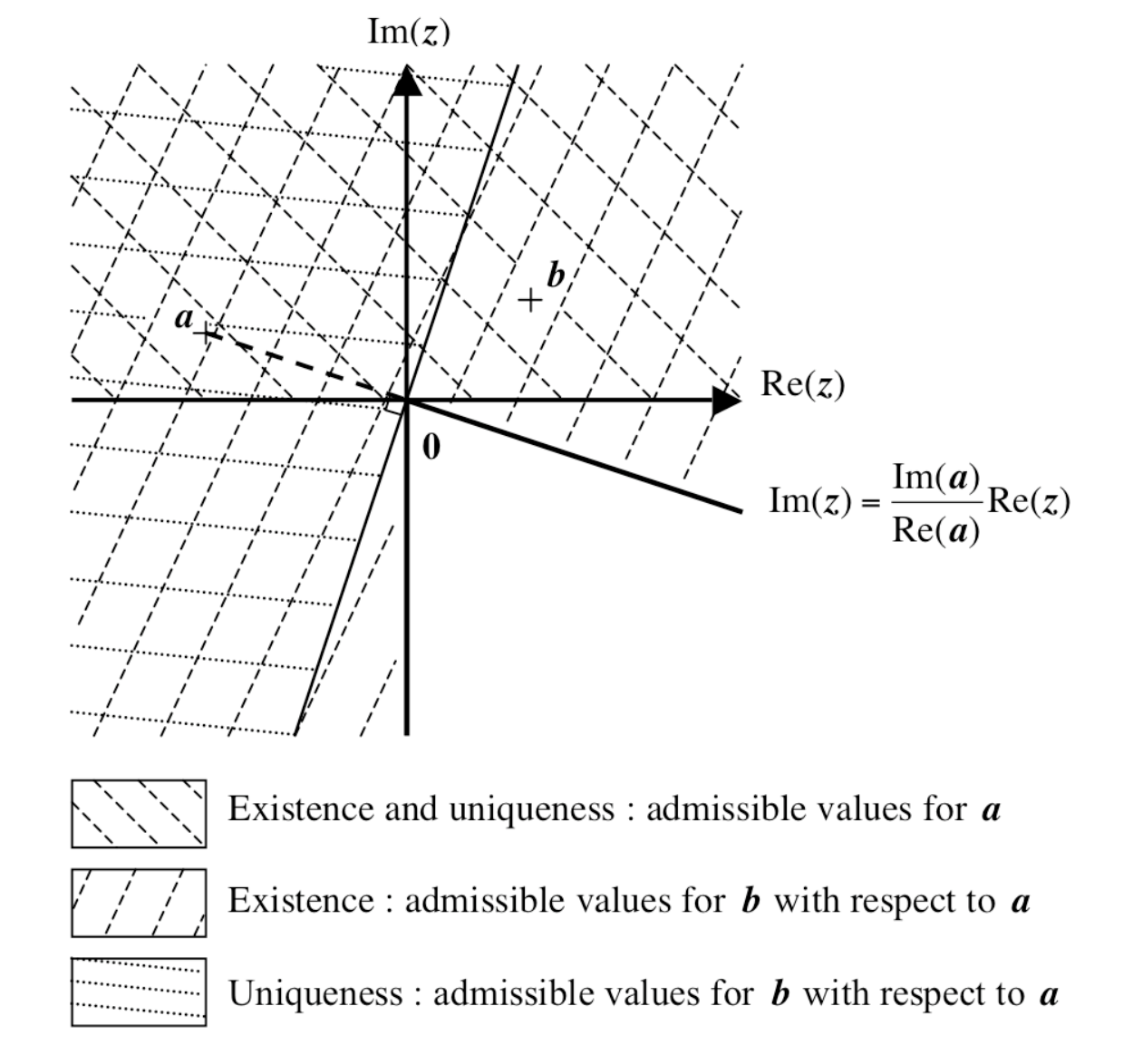}
   \caption{Uniqueness implies existence}
   \label{exiuni}
  \end{center}
 \end{minipage}
\end{figure}

\section{Proofs of the localization properties}
\label{proofsta}

In this Section, we prove Theorems~\ref{thmstag}, \ref{thmF}, \ref{thmstacom}, \ref{bound}, \ref{thmstaO} and \ref{thmcompRN}.

\noindent
We recall some useful Gagliardo-Nirenberg's and Young inequalities.

\begin{prop}
\label{gn}
Let $\Omega\subseteq\R^N$ be a nonempty open subset and let $0\le p\le 1.$ Then, there exists a positive constant $C=C(N)$ such that
\begin{gather}
\label{propgn1}
\forall\u\in\bs{H^1_0}(\Omega)\cap\bs{L^{p+1}}(\Omega), \;\|\u\|_{\bs{L^2}(\Omega)}\le C\|\nabla\u\|_{\bs{L^2}(\Omega)}^\frac{N(1-p)}{(N+2)-p(N-2)}
\|\u\|_{\bs{L^{p+1}}(\Omega)}^\frac{2(1+p)}{(N+2)-p(N-2)}, \medskip \\
\label{propgn2}
\forall\u\in\bs{H^1_0}(\Omega)\cap\bs{L^1}(\Omega), \;
\|\u\|_{\bs{L^{p+1}}(\Omega)}^{p+1}\le C\|\nabla\u\|_{\bs{L^2}(\Omega)}^\frac{2pN}{N+2}
\|\u\|_{\bs{L^1}(\Omega)}^\frac{(N+2)-p(N-2)}{N+2}.
\end{gather}
Note that $C$ does not depend on $\Omega.$
\end{prop}

\begin{lem}
\label{lemyoung}
For any real $x\ge0,$ $y\ge0,$ $\eps>0$ and $p>1,$ one has
\begin{gather}
\label{young}
xy\le\frac{1}{p^\p}\eps^{p^\p}x^{p^\p}+\frac{1}{p}\eps^{-p}y^p.
\end{gather}
\end{lem}

\begin{lem}
\label{lemAB}
Let $(\va,\vb)\in\C^2$ satisfying~\eqref{ab} and let $C_0,$ $C_1,$ $C_2,$ $C_3$ be four nonnegative real numbers satisfying
\begin{gather}
\label{Im}
\big|C_1+\Im(\va)C_2+\Im(\vb)C_3\big|\le C_0,	\\
\label{Re}
\big|\Re(\va)C_2+\Re(\vb)C_3\big|\le C_0.
\end{gather}
Then one has
\begin{gather}
\label{lemestimAB}
0\le C_1+LC_2\le MC_0,
\end{gather}
where the positive constants $L$ and $M$ are defined by~\eqref{L} and \eqref{M}, respectively.
\end{lem}

\begin{proof*}
We split the proof in 6 cases. Let $\delta>0.$ \\
\textbf{Case 1.} $\Im(\va)>0$ and $\Im(\vb)\ge0.$ \\
Then \eqref{lemestimAB} follows from \eqref{Im}. \\
\textbf{Case 2.} $\Im(\va)=0,$ $\Im(\vb)\ge0$ and $\Re(\va)\Re(\vb)\ge0.$ \\
We compute $\eqref{Im}+\sign(\Re(\va))\eqref{Re}$ and then obtain \eqref{lemestimAB}. \\
\textbf{Case 3.} $\Im(\va)\ge0,$ $\Im(\vb)<0$ and $\Re(\va)\Re(\vb)\ge0.$ \\
We compute $\eqref{Im}+\frac{|\Im(\vb)|}{\Re(\vb)}\eqref{Re}$ and then obtain \eqref{lemestimAB}. \\
\textbf{Case 4.} $\Re(\va)\Re(\vb)<0.$ \\
If $\Im(\vb)=0$ then \eqref{ab} implies $\Im(\va)>0,$ which falls into the scope of Case~1. So we may assume $\Im(\vb)\neq0.$ We compute
$\eqref{Im}-\frac{\Im(\vb)}{\Re(\vb)}\eqref{Re}$ and then obtain \eqref{lemestimAB}. \\
\textbf{Case 5.} $\Im(\va)<0,$ $\Im(\vb)\ge0$ and $\Re(\va)\Re(\vb)\ge0.$ \\
We compute $\eqref{Im}+\frac{|\Im(\va)|+\delta}{\Re(\va)}\eqref{Re}$ and then obtain \eqref{lemestimAB}. \\
\textbf{Case 6.} $\Im(\va)<0,$ $\Im(\vb)<0$ and $\Re(\va)\Re(\vb)>0.$ \\
We compute $\eqref{Im}+\max\left\{\frac{|\Im(\va)|+\delta}{\Re(\va)},\frac{|\Im(\vb)|}{\Re(\vb)}\right\}\eqref{Re}$ and then obtain \eqref{lemestimAB}.
\\
This ends the proof.
\medskip
\end{proof*}

\begin{vproof}{of Theorems~\ref{thmstag} and~\ref{thmF}.}
In order to establish our result in all cases of \eqref{ab}, we will adopt the proofs of Theorem~2.1 p.12--18 and Theorem~3.2 p.28--30 of Antontsev,
D{\'{\i}}az and Shmarev~\cite{MR2002i:35001}, which has to be adapted. We denote by $\sigma$ the surface measure on a sphere,
$\rho_2=\rho_0,$ if we are concerned by Theorem~\ref{thmstag} and $\rho_2=\rho_1,$ if we are concerned by Theorem~\ref{thmF}. Assume we have either $\rho_2<\dist(x_0,\partial\Omega)$ $\big(\!\!\iff\ovl B(x_0,\rho_2)\subset\Omega\big)$ or
$\rho_2>\dist(x_0,\partial\Omega).$ The remaining case $\rho_2=\dist(x_0,\partial\Omega)$
$\big(\!\!\iff B(x_0,\rho_2)\subset\Omega$ and $\partial\Omega\cap\S(x_0,\rho_2)\not=\emptyset\big),$ will be treated at the end of the proof\footnote{For simplicity, we assume that $\partial\Omega\neq\emptyset.$ Otherwise, we have
$\Omega=\R^N$ and we only have to treat the first case: $\ovl B(x_0,\rho_2)\subset\Omega.$}. If
$\rho_2>\dist(x_0,\partial\Omega),$ we have $\u\in\bs{H^1_0}(\Omega).$ So we may define
$\wt{\u}\in\bs{H^1_0}\big(\Omega\cup B(x_0,\rho_2)\big)$ satisfying $\wt{\u}_{|\Omega}\in\bs{H^1_0}(\Omega),$ by setting $\wt{\u}=\u,$ in $\Omega$ and $\wt{\u}=\0,$ in $\Omega^\c\cap B(x_0,\rho_2).$ Then
$\nabla\wt{\bs  u}=\nabla\u,$ almost everywhere in $\Omega$ and $\nabla\wt{\u}=\0,$ almost everywhere in $\Omega^\c\cap B(x_0,\rho_2).$ Still if $\rho_2>\dist(x_0,\Omega),$ we denote by $\wt\F$ the extension of $\F$ by
$\0$ in $\Omega^\c\cap B(x_0,\rho_2).$ We now proceed with the proof in 7 steps. \\
{\bf Step~1.} Let $L$ and $M$ be the constants defined by~\eqref{L} and \eqref{M}, respectively. For almost every
$\rho\in(0,\rho_2),$
\begin{gather}
\label{proofthmstag1}
\|\nabla\wt{\u}\|_{\bs{L^2}(B(x_0,\rho))}^2+L\|\wt{\u}\|_{\bs{L^{m+1}}(B(x_0,\rho))}^{m+1}
\le MI(\rho)+MJ(\rho),
\end{gather}
where $I(\rho)=\left|\dsp\int_{\S(x_0,\rho)}\wt{\u}\ovl{\nabla\wt{\u}}.\frac{x-x_0}{|x-x_0|}\d\sigma\right|$ and $J(\rho)=\dsp\int_{B(x_0,\rho)}|\wt\F(x)\ovl{\wt\u(x)}|\d x.$ Moreover, $I,J\in L^1(0,\rho_2).$ \\
From Hölder's inequality, the above discussion and Sobolev's embedding,
\begin{gather*}
\|I\|_{L^1(0,\rho_2)}\le\|\wt{\u}\|_{\bs{H^1}(B(x_0,\rho_2))}^2<\infty, \\
\|J\|_{L^1(0,\rho_2)}\le\|\wt\F\|_{\bs{L^\frac{m+1}{m}}(B(x_0,\rho_2))}\|\wt\u\|_{\bs{L^m}(B(x_0,\rho_2))}<\infty.
\end{gather*}
Let $\rho\in(0,\rho_2)$ For any $n\in\N,$
$n>\frac{1}{\rho},$ we define the cutoff function $\psi_n\in W^{1,\infty}(\R)$ by
\begin{gather*}
\forall t\in\R,\; \psi_n(t)=
 \begin{cases}
  \: 1,			& \mbox{ if } |t|\in\left[0,\rho-\frac{1}{n}\right], \medskip \\
  \: n(\rho-|t|),	& \mbox{ if } |t|\in\left(\rho-\frac{1}{n},\rho\right), \medskip \\
  \; 0,			& \mbox{ if } |t|\in[\rho,\infty),
 \end{cases}
\end{gather*}
and we set for almost every $x\in\Omega\cup B(x_0,\rho_2),$ $\bs{\vphi_n}(x)=\psi_n(|x-x_0|)\wt{\u}(x).$ If
$\rho_2<\dist(x_0,\partial\Omega)$ then $\supp\bs{\vphi_n}\subseteq\ovl B(x_0,\rho)\subset\Omega$ and so
$\bs{\vphi_n}\in\bs{H^1_\c}(\Omega).$ If $\rho_2>\dist(x_0,\partial\Omega)$ then
${\bs{\vphi_n}}_{|\Omega}\in\bs{H^1_0}(\Omega)$ and $\supp\bs{\vphi_n}\subseteq\ovl\Omega\cap\ovl B(x_0,\rho).$ It follows from Definition~\ref{defsols} and Remark~\ref{rmkdefsols}, \ref{rmkdefsols2}. and \ref{rmkdefsols3}., that
$\bs\vphi=\vi{\bs{\vphi_n}}_{|\Omega}$ is an admissible test function and so
\begin{multline*}
\Re\vint_{B(x_0,\rho)}\psi_n(|x-x_0|)\left(|\nabla\wt{\u}|^2-\vi\va|\wt{\u}|^{m+1}-\vi\vb|\wt{\u}|^2\right)\d x \\
=-\Re\vint_{B(x_0,\rho)}\psi_n^\p(|x-x_0|)\ovl{\wt{\u}}\nabla\wt{\u}.\frac{x-x_0}{|x-x_0|}\d x
+\Im\vint_{B(x_0,\rho)}\psi_n(|x-x_0|)\wt{\F}\ovl{\wt{\u}}\d x.
\end{multline*}
Introducing the spherical coordinates $(r,\sigma),$ we get
\begin{align*}
	& \; \left|\Re\vint_{B(x_0,\rho)}\psi_n(|x-x_0|)\left(|\nabla\wt{\u}|^2-\vi\va|\wt{\u}|^{m+1}-\vi\vb|\wt{\u}|^2\right)\d x\right|					\\
   =	& \; \left|\Re\left(n\bs\vint_{\rho-\frac{1}{n}}^\rho\left(\:\vint_{\S(x_0,r)}
		\ovl{\wt {\u}}\nabla\wt{\u}.\frac{x-x_0}{|x-x_0|}\d\sigma\right)\d r\right)-\Im\vint_{B(x_0,\rho)}\psi_n(|x-x_0|)\wt{\F}\ovl{\wt{\u}}\d x\right|		\\
  \le	& \; n\bs\vint_{\rho-\frac{1}{n}}^\rho I(r)\d r+\vint_{B(x_0,\rho)}\psi_n(|x-x_0|)|\wt\F(x)\ovl{\wt\u(x)}|\d x.
\end{align*}
We now let $n\nearrow\infty.$ Using the Lebesgue's dominated convergence Theorem and recalling that
$I\in L^1(0,\rho_2),$ we obtain
\begin{gather}
\label{proofthmstag2}
\left|\|\nabla\wt{\u}\|_{\bs{L^2}(B(x_0,\rho))}^2
+\Im(\va)\|\wt{\u}\|_{\bs{L^{m+1}}(B(x_0,\rho))}^{m+1}+\Im(\vb)\|\wt{\u}\|_{\bs{L^2}(B(x_0,\rho))}^2\right|
\le I(\rho)+J(\rho).
\end{gather}
Proceeding as above with $\bs\vphi={\bs{\vphi_n}}_{|\Omega},$ we get
\begin{gather}
\label{proofthmstag3}
\left|\Re(\va)\|\wt{\u}\|_{\bs{L^{m+1}}(B(x_0,\rho))}^{m+1}+\Re(\vb)\|\wt{\u}\|_{\bs{L^2}(B(x_0,\rho))}^2\right|\le I(\rho)+J(\rho).
\end{gather}
Then Step~1 follows from~\eqref{proofthmstag2}, \eqref{proofthmstag3} and Lemma~\ref{lemAB}. \\
Let us recall and introduce some notations. Let $\tau\in\left(\frac{m+1}{2},1\right]$ and let $\rho\in(0,\rho_2).$ We set
\begin{gather*}
\begin{array}{lll}
E(\rho)=\|\nabla\wt{\u}\|_{\bs{L^2}(B(x_0,\rho))}^2, & b(\rho)=\|\wt{\u}\|_{\bs{L^{m+1}}(B(x_0,\rho))}^{m+1},
                                                                               					& \delta=\frac{k}{2(1+m)},					\medskip \\
\theta=\frac{(1+m)+N(1-m)}{k}\in(0,1), & \ell=\frac{1}{\theta(1+m)},	& \gamma(\tau)=\frac{2\tau-(1+m)}{k}\in(0,1),	\medskip \\
\mu(\tau)=\frac{2(1-\tau)}{k}, & \eta(\tau)=\frac{1-m}{1+m}-\gamma(\tau)>0. &
\end{array}
\end{gather*}
{\bf Step~2.} $E\in W^{1,1}(0,\rho_2),$ for a.e.~$\rho\in(0,\rho_2),$
$E^\p(\rho)=\|\nabla\wt{\u}\|_{\bs{L^2}(\S(x_0,\rho))}^2$ and
\begin{multline}
\label{proofthmstag4}
0\le E(\rho)+b(\rho)\le
CL_1ME^\p(\rho)^\frac{1}{2}\left(E(\rho)^\frac{1}{2}
			+\rho^{-\delta}b(\rho)^\frac{1}{m+1}\right)^\theta b(\rho)^\frac{1-\theta}{m+1} \\
+(2L_1M)^\frac{m+1}{m}\|\wt\F\|_{\bs{L^\frac{m+1}{m}}(B(x_0,\rho))}^\frac{m+1}{m},
\end{multline}
where $C=C(N,m)$ and $L_1=\max\left\{1,\frac{1}{L}\right\}.$ \\
We have the identity $E(\rho)=\dsp\int_0^\rho\left(\dsp\int_{\S(x_0,r)}|\nabla\wt{\u}|^2\d\sigma\right)\d r.$ Since the mapping
$r\longmapsto\dsp\int_{\S(x_0,r)}|\nabla\wt{\u}|^2\d\sigma$ lies in $L^1(0,\rho_2),$ $E$ is absolutely continuous on $(0,\rho_2).$ We then get the first part of the claim and we only have to establish \eqref{proofthmstag4}. Let $\rho\in(0,\rho_2).$ It follows from Cauchy-Schwarz's inequality that
\begin{gather}
\label{proofthmstag5}
I(\rho)\le\|\nabla\wt{\u}\|_{\bs{L^2}(\S(x_0,\rho))}\|\wt{\u}\|_{\bs{L^2}(\S(x_0,\rho))}
=E^\p(\rho)^\frac{1}{2}\|\wt{\u}\|_{\bs{L^2}(\S(x_0,\rho))}.
\end{gather}
We recall the interpolation-trace inequality (see Corollary~2.1 in D\'iaz and Véron~\cite{MR792828}. Note there is a misprint: $\delta$ has to be replaced with $-\delta).$
\begin{gather}
\label{proofthmstag6}
\|\wt{\u}\|_{\bs{L^2}(\S(x_0,\rho))}\le C\left(\|\nabla\wt{\u}\|_{\bs{L^2}(B(x_0,\rho))}
+\rho^{-\delta}\|\wt{\u}\|_{\bs{L^{m+1}}(B(x_0,\rho))}\right)^\theta\|\wt{\u}\|_{\bs{L^{m+1}}(B(x_0,\rho))}^{1-\theta},
\end{gather}
where $C=C(N,m).$ Putting together \eqref{proofthmstag1}, \eqref{proofthmstag5} and \eqref{proofthmstag6}, we obtain,
\begin{gather}
\label{proofthm2F1}
E(\rho)+b(\rho)\le CL_1ME^\p(\rho)^\frac{1}{2}\left(E(\rho)^\frac{1}{2}
+\rho^{-\delta}b(\rho)^\frac{1}{m+1}\right)^\theta b(\rho)^\frac{1-\theta}{m+1}
+L_1M\vint_{B(x_0,\rho)}|\wt\F(x)\ovl{\wt\u(x)}|\d x.
\end{gather}
Applying Young's inequality (Lemma~\ref{lemyoung}) with $x=\|\wt\F\|_{\bs{L^\frac{m+1}{m}}(B(x_0,\rho))},$ $y=\|\wt\u\|_{\bs{L^{m+1}}(B(x_0,\rho))},$
\linebreak
$\eps=\left(\frac{2L_1M}{m+1}\right)^\frac{1}{m+1}$ and $p=m+1,$ we get
\begin{gather}
\label{proofthm2F}
\vint_{B(x_0,\rho)}|\wt\F(x)\ovl{\wt\u(x)}|\d x\le\frac{m}{m+1}\left(\frac{2L_1M}{m+1}\right)^\frac{1}{m}
\|\wt\F\|_{\bs{L^\frac{m+1}{m}}(B(x_0,\rho))}^\frac{m+1}{m}+\frac{1}{2L_1M}b(\rho),
\end{gather}
for any $\rho\in(0,\rho_2).$ Putting together \eqref{proofthm2F1} and \eqref{proofthm2F}, we obtain \eqref{proofthmstag4}. Hence Step~2. \\
{\bf Step~3.} Let $C_0$ be the constant in~\eqref{proofthmstag4}. For any $\tau\in\left(\frac{m+1}{2},1\right]$ and for a.e.~$\rho\in(0,\rho_2),$
\begin{multline}
\label{proofthmstag7}
C_0L_1ME^\p(\rho)^\frac{1}{2}\left(E(\rho)^\frac{1}{2}+\rho^{-\delta}b(\rho)^\frac{1}{m+1}\right)^\theta
															b(\rho)^\frac{1-\theta}{m+1} \\
\le\left(K_1(\tau)\rho^{-(\nu-1)}E^\p(\rho)\right)^\frac{1}{2}\left(E(\rho)+b(\rho)\right)^\frac{\gamma(\tau)+1}{2},
\end{multline}
where $K_1(\tau)=CL_1^2M^2\max\left\{\rho_2^{\nu-1},1\right\}\max\{b(\rho_2)^{\mu(\tau)},b(\rho_2)^{\eta(\tau)}\}$ and $C=C(N,m).$ \\
Let $\tau\in\left(\frac{m+1}{2},1\right]$ and let $\rho\in(0,\rho_2).$ A straightforward calculation yields
\begin{align*}
     & \; \left(E(\rho)^\frac{1}{2}+\rho^{-\delta}b(\rho)^\frac{1}{m+1}\right)b(\rho)^\frac{1-\theta}{\theta(m+1)} \medskip \\
 =  & \; E(\rho)^\frac{1}{2}b(\rho)^\frac{1-\theta}{\theta(m+1)}+\rho^{-\delta}b(\rho)^\frac{1}{\theta(m+1)} \medskip \\
 =  & \; E(\rho)^\frac{1}{2}b(\rho)^{\tau(1-\theta)\ell}b(\rho)^{(1-\tau)(1-\theta)\ell}
            +\rho^{-\delta}b(\rho)^{\frac{1}{2}+\tau(1-\theta)\ell}b(\rho)^{\ell-\tau(1-\theta)\ell-\frac{1}{2}} \medskip \\
\le & \; 2\rho^{-\delta}\max\left\{\rho_2^\delta,1\right\}K_2(\tau)^\frac{1}{\theta}
            \left(E(\rho)+b(\rho)\right)^{\frac{1}{2}+\tau(1-\theta)\ell},
\end{align*}
where $K_2^2(\tau)=\max\{b(\rho_2)^{\mu(\tau)},b(\rho_2)^{\eta(\tau)}\}.$ Hence~\eqref{proofthmstag7} with $K_1(\tau)=4C_0^2L_1^2M^2K_2^2(\tau)\max\left\{\rho_2^{\nu-1},1\right\}.$ \\
{\bf Step~4.} For any $\tau\in\left(\frac{m+1}{2},1\right]$ and for a.e.~$\rho\in(0,\rho_2),$
\begin{gather}
\label{proofthmstag8}
0\le E(\rho)^{1-\gamma(\tau)}\le K_1(\tau)\rho^{-(\nu-1)}E^\p(\rho)
+(4L_1M)^\frac{(m+1)(1-\gamma(\tau))}{m}\|\wt\F\|_{\bs{L^\frac{m+1}{m}}(B(x_0,\rho))}^\frac{(m+1)(1-\gamma(\tau))}{m}.
\end{gather}
Putting together \eqref{proofthmstag4} and \eqref{proofthmstag7}, and applying again Young's inequality \eqref{young} with $p=\frac{2}{\gamma(\tau)+1},$ $\eps=(\gamma(\tau)+1)^\frac{\gamma(\tau)+1}{2},$
$x=\left(K_1(\tau)\rho^{-(\nu-1)}E^\p(\rho)\right)^\frac{1}{2}$ and
$y=\left(E(\rho)+b(\rho)\right)^\frac{\gamma(\tau)+1}{2},$ we obtain
\begin{align*}
     &  \; E(\rho)+b(\rho) \\
\le & \; \left(K_1(\tau)\rho^{-(\nu-1)}E^\p(\rho)\right)^\frac{1}{2}\left(E(\rho)+b(\rho)\right)^\frac{\gamma(\tau)+1}{2}
		+(2L_1M)^\frac{m+1}{m}\|\wt\F\|_{\bs{L^\frac{m+1}{m}}(B(x_0,\rho))}^\frac{m+1}{m}, \\
\le & \; C\left(K_1(\tau)\rho^{-(\nu-1)}E^\p(\rho)\right)^\frac{1}{1-\gamma(\tau)}+\frac{1}{2}(E(\rho)+b(\rho))
		+(2L_1M)^\frac{m+1}{m}\|\wt\F\|_{\bs{L^\frac{m+1}{m}}(B(x_0,\rho))}^\frac{m+1}{m},
\end{align*}
where $C=\frac{p-1}{p}\eps^\frac{p}{p-1}=C(N,m).$ Changing, if needed, the constant $C$ in the definition of $K_1(\tau),$ we obtain
\begin{gather*}
E(\rho)+b(\rho)\le\left(K_1(\tau)\rho^{-(\nu-1)}E^\p(\rho)\right)^\frac{1}{1-\gamma(\tau)}
+(4L_1M)^\frac{m+1}{m}\|\wt\F\|_{\bs{L^\frac{m+1}{m}}(B(x_0,\rho))}^\frac{m+1}{m}.
\end{gather*}
Raising both sides of the above inequality to the power $1-\gamma(\tau)$ and recalling that
$\big(1-\gamma(\tau)\big)\in(0,1),$ we obtain \eqref{proofthmstag8}. \\
{\bf Step~5.} Let $\alpha\in(0,\rho_0].$ If $E(\alpha)=0$ then $\u_{|B_\Omega(x_0,\alpha)}\equiv\0.$ \\
From our hypothesis, $E^\p=0$ on $(0,\alpha).$ Furthermore, $\|\wt\F\|_{\bs{L^\frac{m+1}{m}}(B(x_0,\alpha))}=0$ \big(from assumption of Theorem~\ref{thmstag} or~\eqref{thmF1}\big). It follows from Step~2 and continuity of $b$ that $b(\alpha)=0.$ Hence Step~5 follows. \\
{\bf Step~6.} Proof of Theorem~\ref{thmstag}. \\
Thus $\rho_2=\rho_0$ and $\|\wt\F\|_{\bs{L^\frac{m+1}{m}}(B(x_0,\rho_0))}=0.$ For any $\tau\in\left(\frac{m+1}{2},1\right],$ set
$r(\tau)^\nu=\left(\rho_0^\nu-\nu\frac{K_1(\tau)E(\rho_0)^{\gamma(\tau)}}{\gamma(\tau)}\right)_+$ and let
$\rho_\M=\vmax_{\tau\in\left(\frac{m+1}{2},1\right]}r(\tau).$ Note that definition of $\rho_\M$ coincides with \eqref{thmstag1}. Let
$\tau\in\left(\frac{m+1}{2},1\right].$ We claim that $E(r(\tau))=0.$ Otherwise, $E(r(\tau))>0$ and so $E>0$ on $[r(\tau),\rho_0).$
From~\eqref{proofthmstag8}, one has (we recall that $\gamma(\tau)-1<0),$
\begin{gather}
\label{proofthmstag9}
\mbox{for a.e.~} \rho\in(r(\tau),\rho_0), \;
K_1(\tau)E^\p(\rho)E(\rho)^{\gamma(\tau)-1}\ge\rho^{\nu-1}.
\end{gather}
We integrate this estimate between $r(\tau)$ and $\rho_0.$ We obtain
\begin{gather*}
\nu\frac{K_1(\tau)}{\gamma(\tau)}\left(E(\rho_0)^{\gamma(\tau)}-E(r(\tau))^{\gamma(\tau)}\right)\ge\rho_0^\nu-r^\nu(\tau).
\end{gather*}
By definition of $r(\tau),$ this gives $E(r(\tau))\le0.$ A contradiction, hence the claim. In particular, $E(\rho_\M)=0.$ It follows from Step~5 that $\u_{|B_\Omega(x_0,\rho_\M)}\equiv0,$ which is the desired result. It remains to treat the case where $\rho_0=\dist(x_0,\partial\Omega).$ We proceed as follows. Let $n\in\N,$ $n\ge\frac{1}{\rho_0}.$ We work on $B\left(x_0,\rho_0-\frac{1}{n}\right)$ instead of $B(x_0,\rho_0)$ and apply the above result. Thus
$\u_{|B\left(x_0,\rho_\M^n\right)}\equiv\0,$ where $\rho_\M^n$ is given by \eqref{thmstag1} with $\rho_0-\frac{1}{n}$ in place of $\rho_0.$ We then let $n\nearrow\infty$ which leads to the result. This finishes the proof of
Theorem~\ref{thmstag}. \\
{\bf Step~7.} Proof of Theorem~\ref{thmF}. \\
We have $\rho_2=\rho_1.$ Let $\gamma=\gamma(1)$ and set for any $\rho\in[0,\rho_1],$
$F(\rho)=(4L_1M)^\frac{(m+1)(1-\gamma)}{m}\|\wt\F\|_{\bs{L^\frac{m+1}{m}}(B(x_0,\rho))}^\frac{(m+1)(1-\gamma)}{m}$ and
$K=K_1(1)\rho_0^{-(\nu-1)}.$ Let $E_\star=\left(\frac{\gamma}{2K}(\rho_1-\rho_0)\right)^\frac{1}{\gamma}$ and
$\eps_\star=\frac{1}{2^{p^\p}(4L_1M)^\frac{m+1}{m}}\left(\frac{\gamma}{2K}\right)^p.$ Note that $p=\frac{1}{\gamma}.$ Assume now
$E(\rho_1)<E_\star.$ Applying Step~4 with $\tau=1,$ one has for a.e.~$\rho\in(\rho_0,\rho_1),$
\begin{gather}
\label{proofthmF9}
-KE^\p(\rho)+E(\rho)^{1-\gamma}\le F(\rho).
\end{gather}
Let define the function $G$ by
\begin{gather}
\label{proofthmF10}
\forall\rho\in[0,\rho_1], \; G(\rho)=\left(\frac{\gamma}{2K}(\rho-\rho_0)_+\right)^\frac{1}{\gamma}.
\end{gather}
Then $G(\rho_1)=E_\star,$ $G\in C^1([0,\rho_1];\R)$ \Big(since $\frac1\gamma>2\Big)$ and $G$ satisfies
\begin{gather}
\label{proofthmF11}
\forall\rho\in[0,\rho_1], \; -KG^\p(\rho)+\frac{1}{2}G(\rho)^{1-\gamma}=0, \\
\label{proofthmF12}
E(\rho_1)<G(\rho_1).
\end{gather}
Finally and recalling that $\gamma=\frac{1}{p},$ from our hypothesis \eqref{thmF1} and \eqref{proofthmF10}, one has
\begin{gather}
\label{proofthmF13}
\forall\rho\in(0,\rho_1), \; F(\rho)\le\frac{1}{2}\left(\frac{\gamma}{2K}(\rho-\rho_0)_+\right)^\frac{1-\gamma}{\gamma}=\frac{1}{2}G(\rho)^{1-\gamma}.
\end{gather}
Putting together \eqref{proofthmF9}, \eqref{proofthmF13} and \eqref{proofthmF11}, one obtains
\begin{gather}
\label{proofthmF14}
-KE^\p(\rho)+E(\rho)^{1-\gamma}\le-KG^\p(\rho)+G(\rho)^{1-\gamma}, \mbox{ for a.e.~} \rho\in(\rho_0,\rho_1).
\end{gather}
Now, we claim that for any $\rho\in[\rho_0,\rho_1),$ $E(\rho)\le G(\rho).$ Indeed, if the claim does not hold, it follows from \eqref{proofthmF12} and continuity of $E$ and $G$ that there exist $\rho_\star\in(\rho_0,\rho_1)$ and
$\delta\in(0,\rho_\star-\rho_0]$ such that
\begin{gather}
\label{proofthmF15}
E(\rho_\star)=G(\rho_\star), \\
\label{proofthmF16}
E(\rho)>G(\rho), \; \forall\rho\in(\rho_\star-\delta,\rho_\star).
\end{gather}
It follows from~\eqref{proofthmF14} and \eqref{proofthmF16} that for a.e.~$\rho\in(\rho_\star-\delta,\rho_\star),$
$G^\p(\rho)<E^\p(\rho).$ But, with \eqref{proofthmF15}, this implies that for any $\rho\in(\rho_\star-\delta,\rho_\star),$ $G(\rho)>E(\rho),$ which contradicts \eqref{proofthmF16}, hence the claim. It follows that $0\le E(\rho_0)\le G(\rho_0)=0.$ We deduce with help of the Step~5 that
$\u_{|B_\Omega(x_0,\rho_0)}\equiv\0,$ which is the desired result. It remains to treat the case where $\rho_1=\dist(x_0,\partial\Omega).$ We proceed as follows. Assume $E(\rho_1)<E_\star.$ Then there exists $\eps>0$ small enough such that $\rho_0<\rho_1-\eps$ and
$E(\rho_1)<E_\star(\eps),$ where $E_\star(\eps)=\left(\frac{\gamma}{2K}(\rho_1-\rho_0-\eps)\right)^\frac{1}{\gamma}.$ Since $\eps_\star$ is a non increasing function of $\rho_1,$ we do not need to change its definition. Estimates \eqref{proofthmF9}--\eqref{proofthmF14} holding with
$\rho_1-\eps$ in place of $\rho_1,$ it follows that $E(\rho_0)=0$ and we finish with the help of Step~5. This ends the proof of Theorem~\ref{thmF}.
\medskip
\end{vproof}

\begin{vproof}{of Theorem~\ref{thmstacom}.}
Let $C_0=C_0(N,m)$ be the constant in estimate \eqref{thmstag1} given by Theorem~\ref{thmstag}. We then choose $C=C_0^{-1}$ in \eqref{thmstacom1} and \eqref{thmstacom2}. Using the notations of Theorem~\ref{thmstag} and its proof, we define for any $\tau\in\left(\frac{m+1}{2},1\right],$
\begin{multline*}
 r(\tau)^\nu=\left((2\rho_0)^\nu-C_0M^2\max\left\{1,\frac{1}{L^2}\right\}\max\left\{(2\rho_0)^{\nu-1},1\right\}\right. \\
 \left.\times\frac{E(2\rho_0)^{\gamma(\tau)}
 \max\{b(2\rho_0)^{\mu(\tau)},b(2\rho_0)^{\eta(\tau)}\}}{2\tau-(1+m)}\right)_+,
\end{multline*}
and recall that $\rho_\M=\vmax_{\tau\in\left(\frac{m+1}{2},1\right]}r(\tau).$
Assume \eqref{thmstacom1} holds. Then $\rho_\M\ge\rho_1(1)\ge\rho_0$ and it follows from \eqref{thmstag1} of Theorem~\ref{thmstag} that $b(\rho_0)=0.$ Now assume \eqref{thmstacom2} holds. Since $E(2\rho_0)\le1,$ $b(2\rho_0)\le1$ and $0<\mu(\tau)<\eta(\tau)<1,$ for any $\tau\in\left(\frac{m+1}{2},1\right),$ it follows from definitions of $\rho_1$ and $\rho_\M,$ that
\begin{gather*}
\rho_\M^\nu\ge\rho_1^\nu(1-s)\ge(2\rho_0)^\nu-C_0M^2\min\{1,L^2\}\frac{\max\{(2\rho_0)^{\nu-1},1\}}{1-m-2s}
b(2\rho_0)^{\mu(1-s)}\ge\rho_0^\nu.
\end{gather*}
By \eqref{thmstag1} of Theorem~\ref{thmstag}, $b(\rho_0)=0.$ This concludes the proof.
\medskip
\end{vproof}

\begin{vproof}{of Theorem~\ref{bound}.}
By Definition~\ref{defsols} and of Remark~\ref{rmkdefsols}, \ref{rmkdefsols3}., we can choose $\bs\vphi=\vi\u$ and
$\bs\vphi=\u$ in \eqref{defsols2}. We then obtain,
\begin{gather*}
\|\nabla\u\|_{\bs{L^2}(\Omega)}^2+\Im(\va)\|\u\|_{\bs{L^{m+1}}(\Omega)}^{m+1}
+\Im(\vb)\|\u\|_{\bs{L^2}(\Omega)}^2=\Im\vint_\Omega\F\ovl\u\d x, \medskip \\
\Re(\va)\|\u\|_{\bs{L^{m+1}}(\Omega)}^{m+1}+\Re(\vb)\|\u\|_{\bs{L^2}(\Omega)}^2
=\Re\vint_\Omega\F\ovl\u\d x.
\end{gather*}
Applying Lemma~\ref{lemAB}, these estimates yield,
\begin{gather}
\label{bound3}
\|\nabla\u\|_{\bs{L^2}(\Omega)}^2+L\|\u\|_{\bs{L^{m+1}}(\Omega)}^{m+1}\le M\vint_\Omega\left|\F\right|\left|\u\right|\d x.
\end{gather}
We apply Young's inequality \eqref{young} with $x=|\F|,$ $y=|\u|,$ $\eps=\left(\frac{2M}{(m+1)L}\right)^\frac{1}{m+1}$ and $p=m+1.$ With \eqref{bound3}, we get
\begin{gather*}
\|\nabla\u\|_{\bs{L^2}(\Omega)}^2+\frac{L}{2}\|\u\|_{\bs{L^{m+1}}(\Omega)}^{m+1}
<M\left(\frac{2M}{L}\right)^\frac{1}{m}\|\F\|_{\bs{L^\frac{m+1}{m}}(\Omega)}^\frac{m+1}{m},
\end{gather*}
from which we deduce \eqref{aprioribound1}. Finally, applying Gagliardo-Nirenberg's inequality~\eqref{propgn1}, with $p=m,$ and Young's inequality \eqref{young}, with $p=\frac{4+N(1-m)}{N(1-m)}$ and $\eps=1,$ one obtains
\begin{gather*}
\|\u\|_{\bs{L^2}(\Omega)}^{2\frac{(N+2)-m(N-2)}{4+N(1-m)}}\le C\|\nabla\u\|_{\bs{L^2}(\Omega)}^\frac{2N(1-m)}{4+N(1-m)}
\|\u\|_{\bs{L^{m+1}}(\Omega)}^\frac{4(1+m)}{4+N(1-m)}<C\left(\|\nabla\u\|_{\bs{L^2}(\Omega)}^2+\|\u\|_{\bs{L^{m+1}}(\Omega)}^{m+1}\right),
\end{gather*}
and finally
\begin{gather}
\label{bound4}
\|\u\|_{\bs{L^2}(\Omega)}^2<C\left(\|\nabla\u\|_{\bs{L^2}(\Omega)}^2+\|\u\|_{\bs{L^{m+1}}(\Omega)}^{m+1}\right)^{\delta+1},
\end{gather}
where $\delta=\frac{2(1-m)}{(N+2)-m(N-2)}.$ Estimate \eqref{aprioribound2} then follows from \eqref{aprioribound1} and \eqref{bound4}.
\medskip
\end{vproof}

\begin{vproof}{of Theorem~\ref{thmstaO}.}
Let $C$ be the constant given by Theorem~\ref{thmstacom} and let $\eps>0.$ Set $K=\supp F$ and
$K(\eps)=\ovl{\vO(\eps)}.$ We would like to apply Theorem~\ref{thmstacom} with $\rho_0=\frac{\eps}{4}.$ By \eqref{aprioribound1} of Theorem~\ref{bound}, there exists $\delta_0=\delta_0(\eps,N,m,L,M)>0$ such that if
$\|\F\|_{\bs{L^\frac{m+1}{m}}(\Omega)}\le\delta_0$ then $\|\u\|_{\bs{L^{m+1}}(\Omega)}\le1$ and
\begin{gather}
 \label{proofthmO1}
  \|\nabla\u\|_{\bs{L^2}(\Omega)}^\frac{2(1-m)}{k}
  \le C2^{-2\nu}(2^\nu-1)(1-m)M^{-2}\min\{1,L^2\}\min\{2,\eps\}^{\nu-1}\eps.
\end{gather}
We recall that the distance between two closed sets $\vA$ and $\vB$ of $\R^N$ with one of them compact is defined by
\begin{gather*}
 \dist(\vA,\vB)=\min_{(x,y)\in\vA\times\vB}|x-y|
\end{gather*}
and that
\begin{gather*}
 \dist(\vA,\vB)>0 \iff \vA\cap\vB=\emptyset.
\end{gather*}
Let $x_0\in\ovl{K(\eps)^\c}.$ Let $y\in\ovl B\left(x_0,\frac{\eps}{2}\right)$ and let $z\in K.$ By definition of $K(\eps),$ $\dist(\ovl{K(\eps)^\c},K)=\eps.$ We then have
\begin{gather*}
 \eps=\dist(\ovl{K(\eps)^\c},K)\le|x_0-z|\le|x_0-y|+|y-z|\le\frac{\eps}{2}+|y-z|.
\end{gather*}
Taking the minimum on $(y,z)\in\ovl B\left(x_0,\frac{\eps}{2}\right)\times K,$ we get
\begin{gather*}
 \frac{\eps}{2}\le\dist\left(\ovl B\left(x_0,\frac{\eps}{2}\right),K\right),
\end{gather*}
which means that $\ovl B\left(x_0,\frac{\eps}{2}\right)\cap K=\emptyset,$ for any $x_0\in\ovl{K(\eps)^\c}.$ By \eqref{proofthmO1}, $\u$ satisfies \eqref{thmstacom1} with $\rho_0=\frac{\eps}{4}$ and we deduce that for any $x_0\in\ovl{K(\eps)^\c},$ $\u_{\left|\Omega\cap B\left(x_0,\frac{\eps}{4}\right)\right.}\equiv\0$ (Theorem~\ref{thmstacom}). Let $n\in\N.$ By compactness,
$\ovl{K\left(\frac{7\eps}{8}\right)^\c}\cap\ovl B(0,n)$ may be covered by a finite number of balls
$B\left(x_0,\frac{\eps}{4}\right)$ with $x_0\in\ovl{K(\eps)^\c}.$ Thus for any $n\in\N,$
$\u_{\left|\Omega\cap K\left(\frac{7\eps}{8}\right)^\c\cap B(0,n)\right.}\equiv\0.$ It follows that $\u=\0$ almost everywhere on
\begin{gather*}
 \bigcup_{n\in\N}\left(\Omega\cap K\left(\frac{7\eps}{8}\right)^\c\cap B(0,n)\right)
 =\Omega\cap K\left(\frac{7\eps}{8}\right)^\c.
\end{gather*}
This means that $\supp\u\subset\ovl\Omega\cap K\left(\frac{7\eps}{8}\right)\subset\ovl\Omega\cap\vO(\eps).$ Finally, since $K$ is a compact set, $\Omega$ is open and $K\subset\Omega$, it follows that if $\eps$ is small enough then
$\vO(\eps)\subset\Omega.$ This ends the proof.
\medskip
\end{vproof}

\begin{vproof}{of Theorem~\ref{thmcompRN}.}
Let $L,$ $M$ and $C$ be the constants given by \eqref{L}, \eqref{M} and Theorem~\ref{thmstacom}, respectively. We would like to apply Theorem~\ref{thmstacom} with $\rho_0=1.$ Since $\F$ is compactly supported and $\u\in\bs{H^1}(\R^N)\cap\bs{L^\frac{m+1}{m}}(\R^N)$, there exists $R>1$ such that $\supp\F\subset B(0,R-1),$
\begin{gather*}
\|\u\|_{\bs{L^{m+1}}(\{|x|>R-1\})}\le1 \; \mbox{ and } \; 
\|\nabla\u\|_{\bs{L^2}(\{|x|>R-1\})}^\frac{2(1-m)}{k}\le C2^{1-\nu}(2^\nu-1)(1-m)M^{-2}\min\{1,L^2\}.
\end{gather*}
Let $x_0\in\R^N$ be such that $|x_0|\ge R+1.$ Then
$\ovl B(x_0,2)\cap\supp\F=\emptyset$ and, with help of the above estimate, $\u$ satisfies \eqref{thmstacom1} with $\rho_0=1.$ It follows from Theorem~\ref{thmstacom} that
$\u_{|B(x_0,1)}\equiv\0.$ For each integer $n\ge2,$ define the compact set $C_n$ by
\begin{gather*}
C_n=\left\{x\in\R^N;\;R+\frac{1}{n}\le |x|\le R+n-\frac{1}{n}\right\}.
\end{gather*}
By compactness, $C_n$ may be covered by a finite number of balls $B(x_0,1),$ where $R+1\le|x_0|\le R+1+n.$ Thus for any $n\in\N,$ $\u_{|C_n}\equiv\0.$ It follows that $\u=\0$ almost everywhere on
\begin{gather*}
\bigcup_{n\ge2}C_n=\Big\{x\in\R^N;\;|x|>R\Big\}.
\end{gather*}
Then $\supp\u\subset\ovl B(0,R),$ which is the desired result.
\medskip
\end{vproof}

\section{Proofs of the existence and smoothness results}
\label{proofexi}

In this Section, we prove Proposition~\ref{propregedp}, Theorem~\ref{thmexista} and \ref{corevo}.

\begin{vproof}{of Proposition~\ref{propregedp}.}
By Remarks~\ref{rmkpropedp1}, equation~\eqref{F} makes senses in $\bs{L_\loc^1}(\Omega).$ \\
\textbf{Proof of Property~1).}
Let $1<q\le p<\infty.$ Assume $\F\in\bs{L_\loc^p}(\Omega)$ and $\u\in\bs{L_\loc^q}(\Omega)$ is a solution to~\eqref{F}. For $r\in(1,\infty),$
$r^-$ denotes any real in $(1,r).$ Assume $\v\in\bs{L_\loc^{r^-}}(\Omega),$ for some $1<r<\infty,$ is a solution of~\eqref{F}. It follows that
$|\v|^{-(1-m)}\v\in\bs{L_\loc^\frac{r^-}{m}}(\Omega)$ and since $0<m<1,$ $\bs{L_\loc^\frac{r^-}{m}}(\Omega)\subset\bs{L_\loc^r}(\Omega).$ So by \eqref{F} and Hölder's inequality, $\bs{Vu}\in\bs{L_\loc^{r^-}}(\Omega)$ and so $\Delta\v\in\bs{L_\loc^{\min\{r^-,p\}}}(\Omega).$ Furthermore, if for some $1<r<\infty,$ $\v\in\bs{L_\loc^r}(\Omega;\C)$ and $\Delta\v\in\bs{L_\loc^r}(\Omega;\C)$ then $\v\in\bs{W_\loc^{2,r}}(\Omega;\C)$ (see for instance Cazenave~\cite{caz-sle}, Theorem~4.1.2 p.101--102). We then have shown the following property. Let $1<r<\infty.$
\begin{gather}
 \label{propregedp1}
  \u\in\bs{L_\loc^{r^-}}(\Omega)\implies\u\in\bs{W_\loc^{2,\min\{r^-,p\}}}(\Omega).
\end{gather}
Now, we proceed to the proof of Property~1) in 2 cases. \\
\textbf{Case 1.} $\left(\frac{N}{2}\le q\le p\right)$ or $\left(q<\frac{N}{2}\text{ and }q\le p\le\frac{Nq}{N-2q}\right).$ \\
It follows from~\eqref{propregedp1}, applied with $r=q,$ that $\u\in\bs{W_\loc^{2,q^-}}(\Omega).$ In one hand, if $q<\frac{N}{2}$ then
$\bs{W_\loc^{2,q^-}}(\Omega)\subset\bs{L^{p^-}_\loc}(\Omega).$ It follows from~\eqref{propregedp1} (applied with $r=p)$ and Sobolev's embedding that $\u\in\bs{L^{p+\delta}_\loc}(\Omega),$ for $\delta\in(0,1)$ small enough. On the other hand, if $q\ge\frac{N}{2}$ then
$\bs{W_\loc^{2,q^-}}(\Omega)\subset\bs{L^{p+1}_\loc}(\Omega).$ So in both cases, $\u\in\bs{L^{p+\delta}_\loc}(\Omega).$
Applying~\eqref{propregedp1} with $r=p+\delta,$ we then obtain $\u\in\bs{W_\loc^{2,p}}(\Omega).$ \\
\textbf{Case 2.} $1<q<p,$ $q<\frac{N}{2}$ and $\frac{Nq}{N-2q}<p.$ \\
We recall that if $1<r<\frac{N}{2}$ then Sobolev's embedding is
\begin{gather}
 \label{propregedp2}
  \bs{W_\loc^{2,r^-}}(\Omega)\subset\bs{L_\loc^{s^-}}(\Omega), \text{ for any } \: 1\le s<\infty \text{ such that } \:
  \dfrac{1}{s}\ge\dfrac{1}{r}-\dfrac{2}{N}.
\end{gather}
Since $\frac{Nq}{N-2q}<p,$ we may define the smallest integer $n_0\ge2$ such that $\frac{1}{q}-\frac{2n_0}{N}<\frac1p.$ We then set
\begin{gather*}
\frac{1}{p_{n_0}}=
	\begin{cases}
		\frac{1}{p+1},				&	\text{if } \frac{1}{q}-\frac{2n_0}{N}\le0,	\medskip \\
		\frac{1}{q}-\frac{2n_0}{N},		&	\text{if } \frac{1}{q}-\frac{2n_0}{N}>0,
	\end{cases}
\end{gather*}
in order to have $p<p_{n_0}<\infty.$ Finally, define the $n_0$ real $(p_n)_{n\in\li0,n_0-1\ri}$ by $p_0=q$ and
\begin{gather*}
 \forall n\in\li0,n_0-1\ri, \; \frac{1}{p_n}=\frac{1}{p_0}-\frac{2n}{N}.
\end{gather*}
It follows that for any $n\in\li1,n_0-1\ri,$ $q\le p_{n-1}<p_n\le p<p_{n_0}<\infty$ and
\begin{gather}
 \label{propregedp3}
  \forall n\in\li1,n_0\ri, \; \frac{1}{p_n} \ge \frac{1}{p_{n-1}} - \frac{2}{N}. 
\end{gather}
From~\eqref{propregedp1}--\eqref{propregedp3} applied $n_0$ times (and recalling that $p<p_{n_0}<\infty),$ we then obtain $\u\in\bs{W_\loc^{2,p}}(\Omega).$ This ends the proof of Property~1). \\
\textbf{Proof of Property~2).}
We recall the following Sobolev's embedding and estimate.
\begin{gather}
 \label{propregedp4}
  \bs{W_\loc^{2,N+1}}(\Omega)\subset\bs{C^{1,\frac{1}{N+1}}_\loc}(\Omega)\subset\bs{C^{0,1}_\loc}(\Omega),
  \medskip \\
 \label{propregedp5}
  \forall (\bs{z_1},\bs{z_2})\in\C^2, \;
  \left||\bs{z_1}|^{-(1-m)}\bs{z_1}-|\bs{z_2}|^{-(1-m)}\bs{z_2}\right|\le5|\bs{z_1}-\bs{z_2}|^m.
\end{gather}
Assume further that $(\F,\bs V)\in\bs{C_\loc^{0,\alpha}}(\Omega)\times\bs{C_\loc^{0,\alpha}}(\Omega),$ for some
$\alpha\in(0,m].$ In particular, $\bs V\in\bs{L_\loc^\infty}(\Omega)$ and by Property~1), $\u\in\bs{W_\loc^{2,N+1}}(\Omega).$ It follows from \eqref{propregedp4}, \eqref{propregedp5} and \eqref{F} that $|\u|^{-(1-m)}\u\in\bs{C_\loc^{0,m}}(\Omega)$ and so $\Delta\u\in\bs{C_\loc^{0,\alpha}}(\Omega).$ Thus $\u\in\bs{C^{2,\alpha}_\loc}(\Omega)$ (Theorem~9.19 p.243--244 in Gilbarg and Trudinger~\cite{MR1814364}). This concludes the proof of the proposition.
\medskip
\end{vproof}

\begin{vproof}{of Theorem~\ref{thmexista}.}
Let $L$ and $M$ be the constants given by \eqref{L} and \eqref{M}, respectively. We proceed in 4 steps. \\
{\bf Step~1.} Let $\Omega\subset\R^N$ be an open bounded subset and let $\g\in\bs{L^2}(\Omega).$ Then there exists a unique solution $\u\in\bs{H^1_0}(\Omega)$ of
\begin{gather}
 \label{pthmexista1}
  -\Delta\u=\g, \mbox{ in } \bs{L^2}(\Omega).
\end{gather}
Moreover, there exists a positive constant $C=C(|\Omega|,N)$ such that
\begin{gather}
 \label{pthmexista2}
  \left\|(-\Delta)^{-1}\g\right\|_{\bs{H^1_0}(\Omega)}\le C\|\g\|_{\bs{L^2}(\Omega)}, \; \forall\g\in\bs{L^2}(\Omega).
\end{gather}
In particular, the mapping
$(-\Delta)^{-1}:\bs{L^2}(\Omega)\tends\bs{H^1_0}(\Omega)$ is linear continuous. \\
Existence and uniqueness come from Lax-Milgram's Theorem where the bounded  coercive bilinear form $a$ on
$\bs{H^1_0}(\Omega)\times\bs{H^1_0}(\Omega)$ and the bounded linear functional $L$ on $\bs{H^{-1}}(\Omega)$ are defined by
\begin{gather*}
a(\u,\v)=\Re\vint_\Omega\nabla\u(x).\ovl{\nabla\v(x)}\d x \quad \mbox{ and } \quad
\langle L,\v\rangle_{\bs{H^{-1}},\bs{H^1_0}}=\Re\vint_\Omega\v(x)\ovl{\g(x)}\d x,
\end{gather*}
respectively. Note that $a$ is coercive due to Poincaré's inequality. Taking the $\bs{H^{-1}}-\bs{H^1_0}$ duality product of equation
\eqref{pthmexista1} with $\u$ and applying Poincaré's inequality, we obtain estimate \eqref{pthmexista2} and so continuity of $(-\Delta)^{-1}.$ \\
{\bf Step~2.} Let $\Omega\subset\R^N$ be an open bounded subset, let $0<m<1,$ let $(\va,\vb)\in\C^2$ and let
$\F\in\bs{L^2}(\Omega).$ For each $\ell\in\N,$ define $\bs{f_\ell}=\bs{g_\ell}-\vi\F,$ where
\begin{gather}
\forall\v\in\bs{L^2}(\Omega), \; \bs{g_\ell}(\v)=
 \label{pthmexista3}
  \begin{cases}
   \vi\va|\v|^{-(1-m)}\v+\vi\bs{bv},                                   	&	\mbox{ if } |\v|\le\ell, \medskip \\
   \vi\va\ell^m\dfrac{\v}{|\v|}+\vi\vb\ell\dfrac{\v}{|\v|},	&	\mbox{ if } |\v| > \ell.
  \end{cases}
\end{gather}
Then for any $\ell\in\N,$ there exists at least one solution $\bs{u_\ell}\in\bs{H^1_0}(\Omega)$ of
\begin{gather*}
 -\Delta\bs{u_\ell}=\bs{f_\ell}\left(\bs{u_\ell}\right), \mbox{ in } \bs{L^2}(\Omega).
\end{gather*}
It is clear that $(\bs{f_\ell})_{\ell\in\N}\subset\bs{C(L^2(\Omega);L^2(\Omega))}.$ With the help of Step~1 and the continuous and compact embedding $\bs i:\bs{H^1_0}(\Omega)\inj\bs{L^2}(\Omega),$ we may define a continuous and compact sequence of mappings
$(\bs{T_\ell})_{\ell\in\N}$ of $\bs{H^1_0}(\Omega)$ as follows. For any $\ell\in\N,$ set
\begin{gather*}
\begin{array}{rcccccl}
\bs{T_\ell}:\bs{H^1_0}(\Omega)		& \stackrel{\bs i}{\inj}	& \bs{L^2}(\Omega)			& \xrightarrow{\bs{f_\ell}}
							& \bs{L^2}(\Omega)	& \xrightarrow{(-\Delta)^{-1}}	& \bs{H^1_0}(\Omega) \medskip \\
                                                    \v		&       \longmapsto	&                  \bs i(\v)=\v		&       \longmapsto
                                                    		&     \bs{f_\ell}(\v)      	&                \longmapsto    		& (-\Delta)^{-1}(\bs{f_\ell})(\v)
\end{array}
\end{gather*}
Let $\ell\in\N.$ Let $C$ be the constant in \eqref{pthmexista2} and set
$R=C(|\va|+|\vb|+1)\Big(2\ell|\Omega|^\frac{1}{2}+\|\F\|_{\bs{L^2}(\Omega)}\Big).$ Let
$\v\in\bs{H^1_0}(\Omega).$ It follows from \eqref{pthmexista2} that
\begin{align*}
      &    \|\bs{T_\ell}(\v)\|_{\bs{H^1_0}(\Omega)}=\left\|(-\Delta)^{-1}(\bs{f_\ell})(\v)\right\|_{\bs{H^1_0}(\Omega)}
                                                                                                        \le C\|\bs{f_\ell}(\v)\|_{\bs{L^2}(\Omega)} \\
 \le & \; C(|\va|+|\vb|+1)\Big((\ell^m+\ell)|\Omega|^\frac{1}{2}+\|\F\|_{\bs{L^2}(\Omega)}\Big)\le R.
\end{align*}
Hence, $\bs{T_\ell}\left(\bs{H^1_0}(\Omega)\right)\subset\ovl{\bs B}_{\bs{H^1_0}}(0,R),$ where
$\ovl{\bs B}_{\bs{H^1_0}}(0,R)=\left\{\u\in\bs{H^1_0}(\Omega);\|\u\|_{\bs{H^1_0}(\Omega)}\le R\right\}.$ In a nutshell, $\bs{T_\ell}$ is a continuous and compact mapping from $\bs{H^1_0}(\Omega)$ into itself,
$\ovl{\bs B}_{\bs{H^1_0}}(0,R)$ is a closed convex subset of $\bs{H^1_0}(\Omega)$ and
$\bs{T_\ell}\left(\ovl{\bs B}_{\bs{H^1_0}}(0,R)\right)\subset\ovl{\bs B}_{\bs{H^1_0}}(0,R).$ By the Schauder's fixed point Theorem, $\bs{T_\ell}$ admits at least one fixed point $\bs{u_\ell}\in\ovl{\bs B}_{\bs{H^1_0}}(0,R).$ Hence Step~2 follows. \\
{\bf Step~3.} Let be the hypotheses of the theorem. Assume further that $\Omega$ is bounded. Then
equation~\eqref{snls} admits at least one solution $\u\in\bs{H^1_0}(\Omega).$ \\
In other words, we have to solve
\begin{gather}
 \label{pthmexista4}
  -\Delta\u=\f(\u), \mbox{ in } \bs{L^2}(\Omega),
\end{gather}
where $\f=\g-\vi\F$ and for any $\v\in\bs{L^2}(\Omega),$ $\g(\v)=\vi\va|\v|^{-(1-m)}\v+\vi\bs{bv}.$ Let
$(\bs{F^k})_{k\in\N}\subset\bs\Dr(\Omega)$ be such that
$\bs{F^k}\xrightarrow[k\to\infty]{\bs{L^\frac{m+1}{m}}(\Omega)}\F$ and for any $k\in\N,$
$\|\bs{F^k}\|_{\bs{L^\frac{m+1}{m}}(\Omega)}\le2\|\F\|_{\bs{L^\frac{m+1}{m}}(\Omega)}.$ Let
$\bs{g_\ell}$ be defined by \eqref{pthmexista3} and set for any $(k,\ell)\in\N^2,$ $\bs{f_\ell^k}=\bs{g_\ell}-\vi\bs{F^k}.$ For any
$(k,\ell)\in\N^2,$ let $\bs{u_\ell^k}\in\bs{H^1_0}(\Omega)$ be a solution of
\begin{gather}
 \label{pthmexista5}
  -\Delta\bs{u_\ell^k}=\bs{f_\ell}(\bs{u_\ell^k}), \mbox{ in } \bs{L^2}(\Omega),
\end{gather}
given by Step~2. We take the $\bs{H^{-1}}-\bs{H^1_0}$ duality product of equation \eqref{pthmexista5} with
$\bs{u_\ell^k}$ first and $\vi\bs{u_\ell^k}$ second. Applying Lemma~\ref{lemAB}, we then get for any $(k,\ell)\in\N^2,$
\begin{multline*}
 \|\nabla\bs{u_\ell^k}\|_{\bs{L^2}(\Omega)}^2+L\|\bs{u_\ell^k}\|_{\bs{L^{m+1}}(\{|\bs{u_\ell^k}|\le\ell\})}^{m+1}
 									+L\ell^m\|\bs{u_\ell^k}\|_{\bs{L^1}(\{|\bs{u_\ell^k}|>\ell\})}		\medskip \\
 \le M\int_\Omega|\bs{F^k}||\bs{u_\ell^k}
					|\left(\chi_{\left\{|\bs{u_\ell^k}|\le\ell\right\}}+\chi_{\left\{|\bs{u_\ell^k}|>\ell\right\}}\right)\d x.
\end{multline*}
Applying Young's inequality \eqref{young} to the first term on the right-hand side and the Hölder's inequality to the second term of the right-hand side, we arrive to the following estimate.
\begin{multline}
 \label{pthmexista6}
  2\|\nabla\bs{u_\ell^k}\|_{\bs{L^2}(\Omega)}^2+L\|\bs{u_\ell^k}\|_{\bs{L^{m+1}}(\{|\bs{u_\ell^k}|\le\ell\})}^{m+1}
  +2\|\bs{u_\ell^k}\|_{\bs{L^1}(\{|\bs{u_\ell^k}|>\ell\})}\big(L\ell^m-M\|\bs{F^k}\|_{\bs{L^\infty}(\Omega)}\big)	\medskip \\
  \le M\left(\frac{2M}{L}\right)^\frac{1}{m}\|\bs{F^k}\|_{\bs{L^\frac{m+1}{m}}(\Omega)}^\frac{m+1}{m}
	\le C\|\F\|_{\bs{L^\frac{m+1}{m}}(\Omega)}^\frac{m+1}{m}.
\end{multline}
For any $k\in\N,$ there exists $\ell_k\in\N$ large enough such that
$L\ell_k^m-M\|\bs{F^k}\|_{\bs{L^\infty}(\Omega)}\ge1.$ Moreover, $\Omega$ being bounded, we have
$\bs{L^{m+1}}(\Omega)\inj\bs{L^1}(\Omega).$ So $\left(\nabla\bs{u_{\ell_k}^k}\right)_{k\in\N}$ and
$(\bs{u_{\ell_k}^k})_{k\in\N}$ are bounded in $\bs{L^2}(\Omega)$ and $\bs{L^1}(\Omega),$ respectively. It follows from Gagliardo-Nirenberg's inequality \eqref{propgn2} (applied with $p=1),$ that $(\bs{u_{\ell_k}^k})_{k\in\N}$ is also bounded in $\bs{L^2}(\Omega)$ and so in $\bs{H^1_0}(\Omega).$ Finally, by Rellich-Kondrachov's Theorem, there exists a subsequence $(\bs{u_{\vphi(n)}^n})_{n\in\N}$ of $(\bs{u_{\ell_k}^k})_{k\in\N}$ and $h\in L^2(\Omega;\R),$ such that
\begin{align}
\label{pthmexista7}
& \bs{u_{\vphi(n)}^n}\xrightarrow[n\to\infty]{\bs{L^2}(\Omega)}\u, \medskip \\
\label{pthmexista8}
& \bs{u_{\vphi(n)}^n}\xrightarrow[n\to\infty]{\text{a.e.~in }\Omega}\u, \medskip \\
\label{pthmexista9}
& \left|\bs{u_{\vphi(n)}^n}\right|\le h, \mbox{ for any } n\in\N, \mbox{ a.e.~in } \Omega,
\end{align}
By \eqref{pthmexista8} and \eqref{pthmexista9},
\begin{gather*}
\bs{g_{\vphi(n)}}\Big(\bs{u_{\vphi(n)}^n}\Big)\chi_{\big\{|\bs{u_{\vphi(n)}^n}|\le\vphi(n)\big\}}
\xrightarrow[n\to\infty]{\text{a.e.~in }\Omega}\g(\u), \medskip \\
\forall n\in\N, \;
\left|\bs{g_{\vphi(n)}}\Big(\bs{u_{\vphi(n)}^n}\Big)\right|\le C(h^m+h)\in\bs{L^1}(\Omega), \mbox{ a.e.~in } \Omega.
\end{gather*}
It follows from the dominated convergence Theorem that
\begin{gather}
\label{pthmexista10}
\bs{g_{\vphi(n)}}\Big(\bs{u_{\vphi(n)}^n}\Big)\chi_{\big\{|\bs{u_{\vphi(n)}^n}|
\le\vphi(n)\big\}}\xrightarrow[n\to\infty]{\bs{L^1}(\Omega)}\g(\u).
\end{gather}
In addition, by \eqref{pthmexista7} and Hölder's inequality,
\begin{gather}
\label{pthmexista11}
\left\|\bs{g_{\vphi(n)}}\Big(\bs{u_{\vphi(n)}^n}\Big)\chi_{\big\{|\bs{u_{\vphi(n)}^n}|>\vphi(n)\big\}}\right\|_{\bs{L^1}(\Omega)}
\le\frac C{\vphi(n)}\left(\big\|\bs{u_{\vphi(n)}^n}\big\|_{\bs{L^{m+1}}(\Omega)}^{m+1}
+\big\|\bs{u_{\vphi(n)}^n}\big\|_{\bs{L^2}(\Omega)}^2\right)\xrightarrow{n\to\infty}0.
\end{gather}
Putting together \eqref{pthmexista10} and \eqref{pthmexista11}, we obtain
\begin{gather}
\label{pthmexista12}
\bs{g_{\vphi(n)}}\Big(\bs{u_{\vphi(n)}^n}\Big)\xrightarrow[n\to\infty]{\bs{L^1}(\Omega)}\g(\u).
\end{gather}
Since $\bs{F^n}\xrightarrow{n\to\infty}\F$ in $\bs{L^\frac{m+1}{m}}(\Omega)\inj\bs{L^1}(\Omega),$ we deduce with help of \eqref{pthmexista7} and \eqref{pthmexista12} that
\begin{gather}
\label{pthmexista13}
\Delta\bs{u_{\vphi(n)}^n}\xrightarrow[n\to\infty]{\bs{H^{-2}}(\Omega)}\Delta\u, \medskip \\
\label{pthmexista14}
\bs{f_{\vphi(n)}}\Big(\bs{u_{\vphi(n)}^n}\Big)\xrightarrow[n\to\infty]{\bs{L^1}(\Omega)}\f(\u).
\end{gather}
By \eqref{pthmexista5}, we have for any $n\in\N,$
$-\Delta\bs{u_{\vphi(n)}^n}=\bs{f_{\vphi(n)}^n}\Big(\bs{u_{\vphi(n)}^n}\Big),$ in $\bs{L^2}(\Omega).$ Estimates \eqref{pthmexista13} and \eqref{pthmexista14} allow to pass in the limit in this equation in the sense of
$\bs{\Dr^\p}(\Omega).$ This means that $\u\in\bs{H^1_0}(\Omega)$ is a solution of \eqref{pthmexista4} and since
$\f(\u)\in L^2(\Omega),$ equation \eqref{pthmexista4} makes sense in $\bs{L^2}(\Omega).$ \\
{\bf Step~4.} Conclusion. Under the hypotheses of the theorem, equation~\eqref{snls} admits at least one solution $\u\in\bs{H^1_0}(\Omega)\cap\bs{L^{m+1}}(\Omega)$ and Properties 1)--3) of the theorem hold. \\
For any $n\in\N,$ we write $\Omega_n=\Omega\cap B(0,n).$ Let $n_0\in\N$ be large enough to have
$\Omega_{n_0}\not=\emptyset.$ For each $n>n_0,$ let  $\bs{u_n}\in\bs{H^1_0}(\Omega_n)$ be any solution of \eqref{snls} in $\Omega_n$ given by Step~3, with the external source $\bs{F_n}=\F_{|\Omega_n}.$ We define
$\wt{\bs{u_n}}\in\bs{H^1_0}(\Omega)$ by extending $\bs{u_n}$ by $\0$ in $\Omega\cap B(0,n)^\c.$ Then
$\nabla\wt{\bs{u_n}}=\nabla\bs{u_n},$ almost everywhere in $\Omega_n$ and $\nabla\wt{\bs{u_n}}=\0,$ almost everywhere in $\Omega\cap B(0,n)^\c.$ It follows from \eqref{aprioribound2} of Theorem~\ref{bound} that
$(\bs{u_n})_{n\in\N}$ is bounded in $\bs{H^1_0}(\Omega_n)\cap\bs{L^{m+1}}(\Omega_n),$ or equivalently,
$(\wt{\bs{u_n}})_{n\in\N}$ is bounded in $\bs{H^1_0}(\Omega)\cap\bs{L^{m+1}}(\Omega).$ Up to a subsequence, that we still denote by $(\wt{\bs{u_n}})_{n\in\N},$ there exists $\u\in\bs{H^1_0}(\Omega)\cap\bs{L^{m+1}}(\Omega)$ such that $\wt{\bs{u_n}}\weak\u$ in $\bs{H^1_\w}(\Omega),$ as $n\tends\infty,$ and
$\wt{\bs{u_n}}\xrightarrow[n\to\infty]{\bs{L^{m+1}_\loc}(\Omega)}\u.$ Let $\bs\vphi\in\bs\Dr(\Omega).$ Since
$\wt{\bs{u_n}}\xrightarrow[n\to\infty]{\bs{L^{m+1}_\loc}(\Omega)}\u,$ we have
$|\wt{\bs{u_n}}|^{-(1-m)}\wt{\bs{u_n}}\xrightarrow[n\to\infty]{\bs{L^\frac{m+1}{m}_\loc}(\Omega)}|\u|^{-(1-m)}\u,$ and in particular
\begin{gather}
 \label{pthmexista15}
  \lim_{n\to\infty}\langle\va|\wt{\bs{u_n}}|^{-(1-m)}\wt{\bs{u_n}}
  ,\bs\vphi\rangle_{\bs{L^\frac{m+1}{m}}(\Omega),\bs{L^{m+1}}(\Omega)}
  =\langle\va|\u|^{-(1-m)}\u,\bs\vphi\rangle_{\bs{L^\frac{m+1}{m}}(\Omega),\bs{L^{m+1}}(\Omega)}.
\end{gather}
Recalling that $\u\in\bs{H^1_0}(\Omega)$ and $\wt{\bs{u_n}}\weak\u$ in $\bs{H^1_\w}(\Omega),$ as $n\tends\infty,$ we get with help of
\eqref{pthmexista15},
\begin{multline}
 \label{pthmexista16}
  \lim_{n\to\infty}\Big(\langle\vi\nabla\wt{\bs{u_n}},\nabla\bs\vphi\rangle_{\bs{L^2}(\Omega),\bs{L^2}(\Omega)}
 	+\langle\va|\wt{\bs{u_n}}|^{-(1-m)}\wt{\bs{u_n}}
	,\bs\vphi\rangle_{\bs{L^\frac{m+1}{m}}(\Omega),\bs{L^{m+1}}(\Omega)} \medskip \\
  +\langle\vb\wt{\bs{u_n}},\bs\vphi\rangle_{\bs{L^2}(\Omega),\bs{L^2}(\Omega)}\Big)
  	=\langle-\vi\Delta\u+\va|\u|^{-(1-m)}\u+\vb\u\rangle_{\bs{\Dr^\p}(\Omega),\bs\Dr(\Omega)}.
\end{multline}
Let $n_1>n_0$ be large enough to have $\supp\bs\vphi\subset\Omega_{n_1}.$ Using the basic properties of $\wt{\bs{u_n}}$ described as above and the fact $\bs{u_n}$ is a solution of \eqref{snls} in $\Omega_n,$ we obtain for any $n>n_1,$ $\bs\vphi_{|\Omega_n}\in\bs\Dr(\Omega_n)$ and
\begin{gather*}
 \begin{split}
  0 = & \; \langle-\vi\Delta\bs{u_n}+a|\bs{u_n}|^{-(1-m)}\bs{u_n}+\vb\bs{u_n}-\bs{F_n},
               	\bs\vphi_{|\Omega_n}\rangle_{\bs{\Dr^\p}(\Omega_n),\bs\Dr(\Omega_n)} \medskip \\
     = & \; \big\langle\vi\nabla\bs{u_n},
     				\nabla\big(\bs\vphi_{|\Omega_n}\big)\big\rangle_{\bs{L^2}(\Omega_n),\bs{L^2}(\Omega_n)}
     		+\langle\va|\bs{u_n}|^{-(1-m)}\bs{u_n},\bs\vphi_{|\Omega_n}
               	\rangle_{\bs{L^\frac{m+1}{m}}(\Omega_n),\bs{L^{m+1}}(\Omega_n)} \medskip \\
        & \; \qquad +\langle\vb\bs{u_n},\bs\vphi_{|\Omega_n}\rangle_{\bs{L^2}(\Omega_n),\bs{L^2}(\Omega_n)}
               	-\langle\bs{F_n},\bs\vphi_{|\Omega_n}\rangle_{\bs{L^\frac{m+1}{m}}(\Omega_n),
               	\bs{L^{m+1}}(\Omega_n)} \medskip \\
     = & \; \langle\vi\nabla\wt{\bs{u_n}},\nabla\bs\vphi\rangle_{\bs{L^2}(\Omega),\bs{L^2}(\Omega)}
               	+\langle\va|\wt{\bs{u_n}}|^{-(1-m)}\wt{\bs{u_n}},
		\bs\vphi\rangle_{\bs{L^\frac{m+1}{m}}(\Omega),\bs{L^{m+1}}(\Omega)} \medskip \\
        & \; \qquad +\langle\vb\wt{\bs{u_n}},\bs\vphi\rangle_{\bs{L^2}(\Omega),\bs{L^2}(\Omega)}
               -\langle\F,\bs\vphi\rangle_{\bs{L^\frac{m+1}{m}}(\Omega),\bs{L^{m+1}}(\Omega)},
 \end{split}
\end{gather*}
from which we deduce
\begin{gather}
 \begin{split}
  \label{pthmexista17}
      & \; \langle\vi\nabla\wt{\bs{u_n}},\nabla\bs\vphi\rangle_{\bs{L^2}(\Omega),\bs{L^2}(\Omega)}
             +\langle\va|\wt{\bs{u_n}}|^{-(1-m)}\wt{\bs{u_n}},
             \bs\vphi\rangle_{\bs{L^\frac{m+1}{m}}(\Omega),\bs{L^{m+1}}(\Omega)}
             +\langle\vb\wt{\bs{u_n}},\bs\vphi\rangle_{\bs{L^2}(\Omega),\bs{L^2}(\Omega)} \medskip \\
   = & \; \langle\F,\bs\vphi\rangle_{\bs{L^\frac{m+1}{m}}(\Omega),\bs{L^{m+1}}(\Omega)},
 \end{split}
\end{gather}
for any $n>n_1.$ Passing to the limit in \eqref{pthmexista17}, we get with \eqref{pthmexista16},
\begin{gather*}
\langle-\vi\Delta\u+\va|\u|^{-(1-m)}\u+\vb\u,\bs\vphi\rangle_{\bs{\Dr^\p}(\Omega),\bs\Dr(\Omega)}
=\langle\F,\bs\vphi\rangle_{\bs{\Dr^\p}(\Omega),\bs\Dr(\Omega)}, \; \forall\bs\vphi\in\bs\Dr(\Omega),
\end{gather*}
which is the desired result. Properties~1) and 2) follow from Proposition~\ref{propregedp}. Finally, if $\F$ is spherically symmetric then $\u,$ obtained as a limit, is also spherically symmetric. Indeed, we replace all the functional spaces
$\bs E$ with $\bs{E_\rad}$ and we follow the above proof step~by step. For $N=1,$ this includes the case where $\F$ is an even function. Finally, if $\F$ is an odd function, it is sufficient to work with the space
$\bs{E_\mathrm{odd}}=\{\v\in\bs E; \; \v \mbox{ is odd}\}$ in place of $\bs E.$ Hence Property~3).
\medskip
\end{vproof}

\begin{vproof}{of Corollary~\ref{corevo}.}
Let the assumptions of the corollary be satisfied. Let $\va=-\vi\bs\lambda,$ $\vb=\vi b$ and $\G=-\vi\F.$ Then $(\va,\vb)\in\A\times\B$ satisfies
\eqref{ab} and we may apply Theorem~\ref{thmexista} and Theorem~\ref{thmcompRN} to find a solution $\bs\vphi\in\bs{C^{2,m}_\b}(\R^N)$ of
\eqref{snls} compactly supported for such $\va,$ $\vb$ and $\G.$ It follows that $\bs\vphi$ is a solution to \eqref{bs}. A straightforward calculation show that $\u$ defined by \eqref{sw} is a solution to \eqref{nls*}. This ends the proof.
\medskip
\end{vproof}

\section{Proofs of the uniqueness results}
\label{proofuni}

In this Section, we prove Theorems~\ref{thmA}, \ref{thmB}, \ref{thmdepcon} and \ref{thmuni}, and Corollaries~\ref{corthmexiuni}, \ref{corthmuni} and \ref{corfin}. Let $0<m\le1.$ Set for any $\z\in\C,$ $\f(\z)=|\z|^{-(1-m)}\z,$ where it is understood that $\f(\0)=\0.$ The proof of Theorem~\ref{thmdepcon} relies on the two following lemmas.

\begin{lem}
\label{lemmon}
Let $0<m\le1.$ Then there exists a positive constant $C$ such that
\begin{gather*}
\forall(\bs{z_1},\bs{z_2})\in\C^2, \;
\Re\Big(\big(\f(\bs{z_1})-\f(\bs{z_2})\big)(\ovl{\bs{z_1}-\bs{z_2}})\Big)
\ge C\frac{|\bs{z_1}-\bs{z_2}|^2}{(|\bs{z_1}|+|\bs{z_2}|)^{1-m}},
\end{gather*}
as soon as $|\bs{z_1}|+|\bs{z_2}|>0.$
\end{lem}

\begin{proof*}
We denote by $|\:.\:|_2$ the Euclidean norm in $\R^2.$ From Lemma~4.10, p.264 of D{\'{\i}}az~\cite{MR853732}, there exists a positive constant $C$ such that
\begin{gather*}
\left(|X|_2^{-(1-m)}X-|Y|_2^{-(1-m)}Y\right).(X-Y)\ge C\frac{|X-Y|_2^2}{(|X|_2+|Y|_2)^{1-m}},
\end{gather*}
for any $(X,Y)\in\R^2\times\R^2$ satisfying $|X|_2+|Y|_2>0.$ We apply this lemma with
$
X=
 \left(
  \begin{array}{c}
   \Re(\bs{z_1}) \\
   \Im(\bs{z_1})
  \end{array}
 \right)
$
and
$
Y=
 \left(
  \begin{array}{c}
   \Re(\bs{z_2}) \\
   \Im(\bs{z_2})
  \end{array}
 \right).
$
Note that $|X|_2=|\bs{z_1}|,$ $|Y|_2=|\bs{z_2}|$ and $|X-Y|_2=|\bs{z_1}-\bs{z_2}|.$ The result follows from a direct calculation.
\medskip
\end{proof*}

\begin{cor}
\label{cormon}
Let $0<m\le1.$ Then,
\begin{gather*}
\Re\Big(\big(\f(\bs{z_1})-\f(\bs{z_2})\big)(\ovl{\bs{z_1}-\bs{z_2}})\Big)\ge0,
\end{gather*}
for any $(\bs{z_1},\bs{z_2})\in\C^2.$
\end{cor}

\begin{proof*}
The result is clear if $|\bs{z_1}|+|\bs{z_2}|=0.$ Otherwise, apply Lemma~\ref{lemmon}.
\medskip
\end{proof*}

\begin{rmk}
\label{rmkcormon}
Corollary~\ref{cormon} still holds for any $m>0$ and can be directly obtained as follows. The mapping $f$ (considered as a function from $\R^2$ onto $\R^2)$ is the derivative of the convex function
\begin{gather*}
\begin{array}{rccl}
F: & \R^2 &       \tends      & \R								\\
     & (x,y) & \longmapsto & \frac{1}{m+1}(x^2+y^2)^\frac{m+1}{2}.
\end{array}
\end{gather*}
It follows that $f$ is a monotone function (Proposition~5.5 p.25 of Ekeland and Temam~\cite{MR2000j:49001}).
\end{rmk}

\begin{lem}
\label{lemmin}
Let $\Omega\subseteq\R^N$ be an open subset, let $0<m<1,$ let $(\va,\vb)\in\C^2$ satisfying \eqref{abuni} and let $\bs{F_1},\bs{F_2}\in\bs{L^1_\loc}(\Omega)$ be such that
$\bs{F_1}-\bs{F_2}\in\bs{L^2}(\Omega).$ Let $\bs{u_1},\bs{u_2}\in\bs{H^1_0}(\Omega)\cap\bs{L^{m+1}}(\Omega)$ be two solutions of \eqref{F1} and \eqref{F2}, respectively. Then there exists a positive constant $C=C(N,m)$ satisfying the following property. If $\va\neq\0$ then
\begin{multline}
\label{lemmin1}
\Im(\va)\|\nabla\bs{u_1}-\nabla\bs{u_2}\|_{\bs{L^2}}^2
	+C|\va|^2\vint_\omega\frac{|\bs{u_1}(x)-\bs{u_2}(x)|^2}{(|\bs{u_1}(x)|+|\bs{u_2}(x)|)^{1-m}}\d x
		+\Re\left(\va\ovl\vb\right)\|\bs{u_1}-\bs{u_2}\|_{\bs{L^2}}^2 \\
\le\Re\vint_\Omega\ovl\va\big(\bs{F_1}(x)-\bs{F_2}(x)\big)\big(\ovl{\bs{u_1}(x)-\bs{u_2}(x)}\big)\d x,
\end{multline}
where $\omega=\Big\{x\in\Omega;|\bs{u_1}(x)|+|\bs{u_2}(x)|>0\Big\}.$ If $\va=\0$ then
\begin{gather}
\label{lemmin2}
\Re(\vb)\|\bs{u_1}-\bs{u_2}\|_{\bs{L^2}}^2
=\Re\vint_\Omega\big(\bs{F_1}(x)-\bs{F_2}(x)\big)\big(\ovl{\bs{u_1}(x)-\bs{u_2}(x)}\big)\d x, \\
\label{lemmin3}
\|\nabla\bs{u_1}-\nabla\bs{u_2}\|_{\bs{L^2}}^2+\Im(\vb)\|\bs{u_1}-\bs{u_2}\|_{\bs{L^2}}^2
=\Im\vint_\Omega\big(\bs{F_1}(x)-\bs{F_2}(x)\big)\big(\ovl{\bs{u_1}(x)-\bs{u_2}(x)}\big)\d x.
\end{gather}
\end{lem}

\begin{proof*}
Let $\bs{u_1}$ and $\bs{u_2}$ be two solutions of \eqref{snls} and \eqref{dbc} and set $\u=\bs{u_1}-\bs{u_2}$ and $\F=\bs{F_1}-\bs{F_2}.$ Then $\u$ satisfies
\begin{gather}
\label{prooflemmin1}
-\vi\Delta\u+\va\big(\f(\bs{u_1})-\f(\bs{u_2})\big)+\vb\u=\F,
\mbox{ in } \bs{H^{-1}}(\Omega)+\bs{L^\frac{m+1}{m}}(\Omega).
\end{gather}
Assume $\va\neq\0.$ We take the $\bs{H^{-1}}+\bs{L^\frac{m+1}{m}}-\bs{H^1_0}\cap\bs{L^{m+1}}$ duality product of \eqref{prooflemmin1} with
$\va\u.$ We obtain,
\begin{gather}
\label{prooflemmin2}
\Im(\va)\|\nabla\u\|_{\bs{L^2}}^2+
|\va|^2\langle\f(\bs{u_1})-\f(\bs{u_2}),\u\rangle_{\bs{L^\frac{m+1}{m}},\bs{L^{m+1}}}
+\Re\left(\va\ovl\vb\right)\|\u\|_{\bs{L^2}}^2=\langle\ovl\va\F,\u\rangle_{\bs{L^2},\bs{L^2}}.
\end{gather}
Applying Lemma~\ref{lemmon}, there exists a positive constant $C=C(N,m)$ such that
\begin{align}
\label{prooflemmin3}
\langle\f(\bs{u_1})-\f(\bs{u_2}),\u\rangle_{\bs{L^\frac{m+1}{m}},\bs{L^{m+1}}}
\ge C\vint_\omega\frac{|\u(x)|^2}{(|\bs{u_1}(x)|+|\bs{u_2}(x)|)^{1-m}}\d x.
\end{align}
Then \eqref{lemmin1} follows from \eqref{prooflemmin2} and \eqref{prooflemmin3}. We turn out the case $\va=\0.$ Taking the $\bs{H^{-1}}+\bs{L^\frac{m+1}{m}}-\bs{H^1_0}\cap\bs{L^{m+1}}$ duality product of \eqref{prooflemmin1} with $\u$ and $\vi\u,$ one respectively obtains \eqref{lemmin2} and \eqref{lemmin3}.
\medskip
\end{proof*}

\begin{vproof}{of Theorem~\ref{thmdepcon}.}
Note that since $(\va,\vb)\in\C^2\setminus\{(\0,\0)\}$ satisfies \eqref{abuni}, if $\va=\0$ and $\Re(\vb)=0$ then one necessarily has $\Im(\vb)>0.$ We apply estimates \eqref{lemmin1}--\eqref{lemmin3} of Lemma~\ref{lemmin}, according to the different cases, and Cauchy-Schwarz's inequality. Estimates \eqref{thmdepcon1} and \eqref{thmdepcon2} follow.
\medskip
\end{vproof}

\begin{vproof}{of Theorem~\ref{thmuni}.}
Let $\F\in\bs{L^1_\loc}(\Omega)$ and let $\bs{u_1},\bs{u_2}\in\bs{H^1_0}(\Omega)\cap\bs{L^{m+1}}(\Omega)$ be two solutions of~\eqref{snls} and \eqref{dbc}. By Lemma~\ref{lemmin}, \eqref{lemmin1}--\eqref{lemmin3} hold with $\bs{F_1}-\bs{F_2}=\0.$ We first note that, since
$\bs{u_1}-\bs{u_2}\in\bs{H^1_0}(\Omega),$ if $\|\nabla\bs{u_1}-\nabla\bs{u_2}\|_{\bs{L^2}}=0$ then $\bs{u_1}-\bs{u_2}=\0,$ a.e.~in $\Omega$ and uniqueness holds. It follows from hypotheses~\eqref{abuni} and Lemma~\ref{lemmin} that one necessarily has $\|\bs{u_1}-\bs{u_2}\|_{\bs{L^2}}=0,$ $\|\nabla\bs{u_1}-\nabla\bs{u_2}\|_{\bs{L^2}}=0$ or $\vint_\omega\frac{|\bs{u_1}-\bs{u_2}|^2}{(|\bs{u_1}(x)|+|\bs{u_2}(x)|^{1-m})}\d x,$ where
$\omega=\big\{x\in\Omega;|\bs{u_1}(x)|+|\bs{u_2}(x)|>0\big\}.$ Those three cases imply that $\bs{u_1}=\bs{u_2},$ a.e.~in $\Omega.$ This finishes the proof of the theorem.
\medskip
\end{vproof}

\begin{vproof}{of Corollary~\ref{corthmexiuni}.}
Apply Theorem~\ref{thmexista}, Theorem~\ref{thmuni} and Remark~\ref{rmkuniexi}.
\medskip
\end{vproof}

\begin{vproof}{of Corollary~\ref{corthmuni}.}
By uniqueness (Theorem~\ref{thmuni}), $\u\equiv\0$ is the unique solution.
\medskip
\end{vproof}

\begin{vproof}{of Corollary~\ref{corfin}.}
Apply Theorem~\ref{thmcompRN}, Theorem~\ref{thmexista}, Proposition~\ref{propregedp}, Theorem~\ref{thmuni} and Remark~\ref{rmkuniexi}.
\medskip
\end{vproof}

\begin{vproof}{of Theorem~\ref{thmA}.}
Apply Theorem~\ref{thmcompRN} and Corollary~\ref{corthmexiuni}.
\medskip
\end{vproof}

\begin{vproof}{of Theorem~\ref{thmB}.}
Apply Theorem~\ref{thmstaO} and Corollary~\ref{corthmexiuni}.
\medskip
\end{vproof}

\baselineskip .4cm

\bibliographystyle{abbrv}
\bibliography{BiblioPaper7}
\addcontentsline{toc}{section}{References}

\end{document}